\documentclass[opre,nonblindrev]{informs3Custom} 

\OneAndAHalfSpacedXI 

\usepackage{algorithm}
\usepackage{algorithmic}
\usepackage{multirow}
\usepackage{hhline}
\usepackage{appendix}

\usepackage{endnotes}
\usepackage{accents}

\newcommand{\ubar}[1]{\underaccent{\bar}{#1}}

\let\footnote=\endnote
\let\enotesize=\normalsize
\def\notesname{Endnotes}%
\def\makeenmark{$^{\theenmark}$}
\def\enoteformat{\rightskip0pt\leftskip0pt\parindent=1.75em
  \leavevmode\llap{\theenmark.\enskip}}

\usepackage{natbib}
 \bibpunct[, ]{(}{)}{,}{a}{}{,}%
 \def\bibfont{\small}%
 \def\bibsep{\smallskipamount}%
 \def\bibhang{24pt}%
\TheoremsNumberedThrough     
\ECRepeatTheorems

\EquationsNumberedThrough    

\begin{document}


\RUNAUTHOR{Durante, Nascimento, and Powell}

\RUNTITLE{Risk Directed Importance Sampling in SDDP for Grid Level Energy Storage}

\TITLE{Risk Directed Importance Sampling in Stochastic Dual Dynamic Programming with Hidden Markov Models for Grid Level Energy Storage}

\ARTICLEAUTHORS{%
\AUTHOR{Joseph L. Durante}
\AFF{Department of Electrical Engineering, Princeton University, Princeton, NJ 08540, \EMAIL{jdurante@princeton.edu}} 
\AUTHOR{Juliana Nascimento}
\AFF{Department of Operations Research, Princeton University, Princeton, NJ 08540, \EMAIL{jnascime@princeton.edu}}
\AUTHOR{Warren B. Powell}
\AFF{Department of Operations Research, Princeton University, Princeton, NJ 08540, \EMAIL{powell@princeton.edu}}
} 

\ABSTRACT{%
Power systems that need to integrate renewables at a large scale must account for the high levels of uncertainty introduced by these power sources. This can be accomplished with a system of many distributed grid-level storage devices. However, developing a cost-effective and robust control policy in this setting is a challenge due to the high dimensionality of the resource state and the highly volatile stochastic processes involved. We first model the problem using a carefully calibrated power grid model and a specialized hidden Markov stochastic model for wind power which replicates crossing times. We then base our control policy on a variant of stochastic dual dynamic programming, an algorithm well suited for certain high dimensional control problems, that is modified to accommodate hidden Markov uncertainty in the stochastics. However, the algorithm may be impractical to use as it exhibits relatively slow convergence. To accelerate the algorithm, we apply both quadratic regularization and a risk-directed importance sampling technique for sampling the outcome space at each time step in the backward pass of the algorithm. We show that the resulting policies are more robust than those developed using classical SDDP modeling assumptions and algorithms.
}%


\HISTORY{Version as of January 22, 2020}
\maketitle

%


\section{Introduction}
\label{Intro}

Power systems with renewables must account for the high volatility and intermittency of wind and solar energy sources. At relatively low levels of renewable penetration, this uncertainty can often be managed using careful scheduling and ramping of generators along with a backup supply of reserves. At current levels of renewable penetration, this is an effective method to ensure system stability. However, rising renewable penetration levels means that the variability of renewables will have a much larger impact on the system. Specifically, in order to maintain the same level of system reliability at high penetration levels, the amount of spinning reserves must increase dramatically. \cite{simao2017challenge} shows that many gigawatts of spinning reserves must be available (but hardly used) to accommodate a comparable amount of renewable power in the system. Considering that infrastructure must be upgraded to accomplish this as well, this is not a cost-effective strategy.

Alternatively, distributed grid-level storage can assist in smoothing out the variability in renewables. Energy may be stored during periods when there is excess renewable power and released during periods when renewables significantly underperform forecasts to prevent outages. Energy storage on a large scale seems to be a possibility in the near future with improvements in battery storage technology and the introduction of electric vehicles, which may act as storage devices, onto the grid.


The simultaneous optimization of distributed, grid-controlled battery storage and generation ramping decisions in the real time market in a power grid with high-penetrations of renewables is thus an important topic. Other authors have studied the interaction between storage and renewables, but often fail to model the system properly by ignoring system dynamics and only considering the aggregate generation, renewables, demands, and storage capabilities of the system \citep[e.g.][]{Jacobson15060}. This paper presents an efficient and effective algorithm to develop control policies for a realistically modeled system.



The problem of optimizing storage over the grid in real time is a large multiperiod stochastic optimization problem with a high-dimensional resource state. Stochastic dual dynamic prograrmming (SDDP), first developed by \cite{pereira1991multi}, has attracted considerable attention for the type of stochastic linear program that describes our energy storage problem \citep[see][]{gorenstin1992stochastic, shapiro2011analysis, shapiro2013worst, guigues2014sddp, de2015improving, bandarra2017multicut}, with a focus on the generation and management of Benders cuts. Significantly less attention has been given to the modeling of the stochastic processes involved in the problem and its impact on the resulting policy. Specifically, in the context of grid storage problems with high penetrations of renewables, models for wind power forecast errors that do not adequately capture the {\it crossing times} of stochastic processes may lead to the development of less robust policies.

Crossing times are contiguous blocks of time for which a stochastic process is above or below some reference series, such as a forecast. Among the models that do a poor job of capturing crossing time behavior are those that assume intertemporal or interstage independence in any random quantities, which is a property required by classical SDDP \citep{pereira1991multi, philpott2008convergence}. These include independently and identically distributed (IID) models and autoregressive moving average (ARMA) time series models (which can be included in the SDDP formulation under certain circumstances and with an expansion of the resource state). These are popular in the study of energy storage systems \citep[e.g.][]{lohndorf2010optimal, zhou2013managing} and SDDP applications in general \citep[e.g.][]{pereira1989optimal, shapiro2011analysis, lohmann2016spatio}, and both often underestimate crossing times \citep[see][]{durante2017}.

To understand why this property is important, imagine sizing a battery to smooth out variability from a wind power source. Simply replicating the forecast error distribution, but underestimating (or overestimating) the length of the periods for which wind is above or below its forecast, could lead to dramatically underestimating (or overestimating) how large that battery must be. Now consider optimizing a storage system using a wind power model which underestimates the periods of time for which wind will be below its forecast. This will result in control policies that do not maintain sufficient storage levels to handle extended periods when wind is below the forecast.

This paper addresses these concerns by modeling wind power generation with the univariate crossing state hidden semi-Markov model presented in \cite{durante2017} which explicitly incorporates crossing times into the model. It is unique in its ability to capture both crossing time distributions and the distribution of forecast errors. Additionally, the distribution of areas above and below the forecast (the surpluses and deficits of energy produced versus expected output) are accurately replicated by the crossing state model. In \cite{durante2017b} we see that modeling these characteristics is indeed important in storage problems as control policies based on value function approximations (VFAs) for a single storage device developed assuming a crossing state model for wind power are shown to be more robust than those trained with standard ARMA models. Here we look to apply this model to produce robust policies in a high dimensional setting with many storage devices spread across a power grid. 

Incorporating the crossing state stochastic model in our optimization algorithm is not straightforward as there exists interstage dependencies between random quantities under this modeling assumption. Various methods for solving multistage stochastic programming problems relax the assumption of interstage independence \citep[ex.][]{infanger1996cut,lohndorf2013optimizing,lohndorf2015optimal}. We base our solution on the SDDP algorithm from \cite{asamov2015regularized} as it is most appropriate for our problem setting. This version of SDDP allows for general Markov uncertainty in the stochastics by fitting a set of VFAs to each \textit{information state} of the model at each time step. A further modification must be made to the algorithm to allow for use of hidden Markov models. This idea has received some recent attention in \cite{dowsonapartially} as well, which independently developed a variant of SDDP that incorporates belief states for multiperiod stochastic programming using the concept of policy graphs. The work shows, using a partially observable multiperiod inventory management problem as an example, that the performance of SDDP with belief states can be comparable to the performance of an SDDP-based policy which was allowed to have full visibility of the underlying state. The work in this paper differs from \cite{dowsonapartially} in that we consider the time-dependent case and the implementation is different; this is explained further in section \ref{SDDP Algo}.

Due to the nature of the problem, an implementation of unregulated SDDP would converge far too slowly to be used in practice. To accelerate convergence, we include quadratic regularization which penalizes deviations from the incumbent solution in early iterations of the algorithm. This encourages the solution to remain near points in the resource state space that have been previously visited as the VFAs are more accurate in these regions. To further speed up the algorithm, we employ a sampling method to accelerate the backward pass of SDDP. Standard sampling, in which the outcome space is sampled according to the true probability of each outcome occurring, will accelerate numerical convergence, but at some cost to solution quality as risky regions of the sample space may not be sampled often enough. As the central reason behind using the more sophisticated stochastic model is to produce more robust solutions, this is counterproductive. Thus, a risk-directed importance sampling technique for the backward pass of the SDDP algorithm is utilized to achieve a solution that is nearly as robust as the unsampled version, in a fraction of the computation time.

Importance sampling has been incorporated in SDDP algorithms previously, as in \cite{Kozmík2015} where it is utilized in the forward pass to improve the upper bound estimator in a risk-averse setting. Here we look to apply importance sampling in the backward pass to sample outcomes that produce large objective values with higher probability. We develop an adaptive learning algorithm that forms sampling distributions based on observations of outcomes and the resulting value of being in a certain storage state over successive iterations of the algorithm. This draws on work from \cite{jiang2017risk}, which presented the method for a different class of approximate dynamic programming algorithms. Implementing a risk-directed importance sampling scheme within SDDP presents its own unique challenges as we are unlikely to visit the same point in a high-dimensional resource state more than once at each time period over the iterations of the algorithm. We must account for this when learning sampling distributions as the VFAs for the resource state can vary greatly throughout the state space and can also be quite inaccurate during the early iterations of the algorithm.

Much of this paper is focused on the reduction of risk and thus it would naturally follow that our objective function is risk-averse as well. Many approaches within dynamic programming seek to minimize an appropriate time consistent and coherent dynamic risk measure \citep[see]{ruszczynski2010risk}. \cite{philpott2012dynamic}, \cite{guigues2012sampling}, and \cite{shapiro2013risk} all apply this concept in various forms to an SDDP algorithm. However, one practical measure of risk that does not satisfy the time consistency and coherency properties is an end-of-horzion penalty that is only realized if the cumulative sum of some metric exceeds some threshold. The need for this type of measure is observed throughout the engineering fields. Suppose grid operators (and their customers) are willing to tolerate only so many shortages per day and no more. In this case, such a threshold penalty is very useful. In order to incorporate this into our problem, we augment the state variable to track shortages by including this ``shortage state" in our resource state and apply the threshold penalty in the final time period as part of the objective function. This approach -- applying a utility function to additional state variables that track metrics of interest to account for risk -- is seen in other dynamic programming applications such as \cite{bauerle2011markov} and \cite{8264389}, for example. Our application of this method provides the system operator with an intuitive means of controlling exposure to shortages via adjustments of the aforementioned threshold and an associated penalty parameter.


This paper makes the following contributions: 1) We extend SDDP to handle hidden Markov states at each time step, which complicates the backward pass of traditional SDDP. 2) We introduce a risk-directed importance sampling scheme for SDDP and show that it accelerates convergence while producing solutions with lower risk. 3) We show risk can be reduced in storage problems with a high-dimensional resource state by modeling the stochastics with a model that captures crossing times.

The remainder of this paper is organized as follows. Section \ref{Modeling} formally describes the grid-level storage problem by defining the five elements of the stochastic optimization problem. We discuss related SDDP approaches and the development of our proposed solution algorithm in Section \ref{SDDP Algo}, which also includes a detailed discussion of risk-directed importance sampling. This section also has a detailed discussion of risk-directed importance sampling within this application. Numerical results are presented in Section \ref{Numerical Results} and the paper is concluded in Section \ref{Conclusion}.

\section{Mathematical Model of the Energy Storage Problem}
\label{Modeling}

Currently, a three step planning process that works on different time scales is employed by PJM and other U.S. grid operators for scheduling energy generation and controlling transmission:

\begin{henumerate}
\item A day-ahead unit-commitment problem is solved using renewable forecasts for the subsequent day and determines the on/off scheduling of steam generation units, which require relatively large notification times, over the course of the day. This is a large mixed integer program and can take several hours to solve.

\item In the intermediate term, planning occurs every 15 to 30 minutes (depending on the system) to determine which gas turbines will be used. As steam generation unit commitment decisions are already fixed at this step, this is a smaller mixed integer program and can be solved in the smaller time window.

\item Finally, real-time economic dispatch takes place every five minutes. This controls the ramping of generator outputs to ensure the stability of the grid under current system conditions. This large linear program must be solved within each five minute time window.
\end{henumerate} In addition, reserve capacity in the form of either spinning reserve, which can be immediately called upon, or non-spinning reserve, which are typically gas turbines with very fast start-up times, is used to account for any unexpected behavior from renewable sources.

Ideally we would incorporate the distributed storage devices in the unit commitment solutions. However, this would be quite complex as it would require integrating a high dimensional resource variable into an already large mixed integer program. Instead, we use the unit commitment model to determine the integer variables (when steam generators and gas turbines must be on or off), and then incorporate storage into the real time model which can just handle the ramping of generators that are already turned on.  This allows our model to decide whether we should store energy, or adjust the generators. Details on this step are left for Section \ref{Numerical Results}.

The remainder of this section provides a complete model of the stochastic optimization problem for grid storage by defining the state variable, the decision variable, exogenous information, the transition function, and the objective function. We use dynamic programming notation adopted from \cite{powell2011approximate} and \cite{powell2016unified} to describe the problem.

\subsection{The State Variable}
\label{StateVar}

We first define the system state variables. These can be divided into static state variables which do not vary over time, and dynamic state variables which evolve over time. Following the modeling convention from \cite{powell2016unified}, we place all static variables in the initial state variable $S_0$. This contains all data pertinent to the problem, including constants, deterministic variables, and information such as the system configuration. All variables mentioned in the subsequent paragraph belong in the initial state variable.

The power grid forms a graph $\mathbf{G} = (\mathcal{N} , \mathcal{E})$, where the grid buses comprise the set of nodes $\mathcal{N}$ and transmission lines form the set of edges $\mathcal{E}$. Edge $(i, j) \in \mathcal{E}$ represents a transmission line that connects bus $i \in \mathcal{N}$ to bus $j \in \mathcal{N}$. Within the set of nodes $\mathcal{N}$, we have power sources and sinks with different characteristics. Generators each have a power capacity and generation cost. Storage devices are characterized by bounds on energy capacity, charging and discharging efficiency, and their variable storage cost. Finally, loads are simply represented by the demand at time $t$. To describe these, the following notation (which resembles the notation for similar grid storage problems in \cite{asamov2016sddp}) is used:
\begin{hitemize}
\item $\mathcal{G}$ = The set of fossil generators.
\item $\mathcal{B}$ = The set of storage devices.
\item $\mathcal{H}$ = The set of renewable sources.
\item $\kappa_g^l,\kappa_g^u$ = The minimum and maximum power capacities for electricity generator $g \in \mathcal{G}$.
\item $\kappa_b^l,\kappa_b^u$ = The minimum and maximum energy capacities for storage devices $b \in \mathcal{B}$.
\item $\eta_b^+,\eta_b^-$ = The charging and discharging multipliers (efficiencies) for each storage device $b \in \mathcal{B}$.
\item $c_{t,g}^G$ = The variable generation cost in \$/MWh for generator $g \in \mathcal{G}$ at time $t$.
\item $c_{t,b}^B$ = The variable storage cost in \$/MWh for storage devices $b \in \mathcal{B}$ at time $t$.
\item $\ell_{ij}$ = The nominal line capacity in MW for transmission line $(i, j ) \in \mathcal{E}$.
\item $d_t$ = The vector of electricity demands (loads), in MW, for each node at time $t = 0, . . ., T$. To simplify the problem, we assume that electricity demand evolves deterministically according to a demand profile, and thus this belongs in the initial state variable. Demand profiles are taken from actual PJM data collected over the course of 2013 in 5 minute intervals.
\item $f_t^E$ = The vector of renewable forecasts at each time $t =0, ..., T$. For simplicity, we assume that the forecast utilized in the day-ahead unit commitment problem is fixed over the problem horizon. Though in reality these may evolve over time, under this assumption the forecast belongs in the initial state variable as well.
\item $Z_{t,g}^G$ = The boolean variable reflecting the unit commitment decision for generator $g \in \mathcal{G}$ at time $t$. These belong in $S_0$ as they are determined prior to the storage optimization step.
\item $\mathcal{Y}$ = The closed, convex set describing the DC power flow model of the grid following Kirchhoff's laws. It takes into account the structure of the electrical grid and bounds on phase angles.
\end{hitemize} In addition to these static variables, $S_0$ contains the initial state of the dynamic state variables, which are described next.

We now define the dynamic pre-decision state variable, $S_t$, which is a minimally dimensioned function of history necessary to model the system from time $t$ onward. The following variables, which will be further elaborated throughout Section \ref{Modeling}, compose the dynamic pre-decision state variable:

\begin{hitemize}
\item $R_t=(R_t^B,R_t^S)$ = The $|\mathcal{B}| + 1$ dimensional resource state vector comprised of:
\begin{hitemize}
\item $R_t^B$ = $(R_{tb}^{B})_{b\in\mathcal{B}}$, a $|\mathcal{B}|$ dimensional vector containing the energy present in each of the storage devices at time $t$, and 
\item $R_t^S$ = The \textit{shortage state}, which tracks the total shortages (in MWh) observed in the system up through time $t$. This augmentation to the resource state allows us to apply a threshold penalty to the cumulative system shortages at the end of the time horizon (see Section \ref{Objective Function}).
\end{hitemize}
\item $I_t$ = The information state, comprised of:
\begin{hitemize}
\item $E_t^W$ = The current power produced by offshore wind at time $t$, and
\item $I_t^W$ = The information state of the stochastic wind model. The wind model has relatively few discrete information states, and these determine the distribution of the next wind power forecast error (see Section \ref{Exo Info}).
\end{hitemize}
\item $K_t$ = The knowledge state. The information state $I_t^W$ of the stochastic wind model is partially hidden (see Section \ref{Exo Info}). Given the observed values of wind through time $t$, $K_t$ contains our current beliefs about the probability the wind process is in each possible information state.
\end{hitemize}

We also define a dynamic post-decision state variable $S_t^x$, which carries only the information necessary to transition to $S_{t+1}$ \textit{after} a decision has been made \citep{powell2011approximate}. Canonically, this is given by $S_t^x=(R_t^x, I_t^x, K_t^x)$ where $R_t^x$, $I_t^x$, and $K_t^x$ are the post-decision resource, information, and knowledge states respectively. However, in this problem $K_t^x=K_t$ as the knowledge state is not affected by the decision $x_t$, so we can write $S_t^x=(R_t^x, I_t^x, K_t)$. To define $I_t^x$, we note that while $E_t^W$ is needed in the pre-decision state to make a decision, it is does not affect the transition to time $t+1$ under this model (see Sections \ref{Exo Info} and \ref{Transition Function}). Thus, $I_t^x$ is equal to $I_t^W$, the information state of the wind process. Finally, we let $R_t^x=(R_t^{B,x},R_t^{S,x})$ contain the energy in each battery and the cumulative shortages through time $t$ after decision $x_t$ has been made.

\subsection{The Decision Variables}
The decision variables at each time $t$ are given by
\begin{align}
x_t=(x_t^{G+},x_t^{B+},x_t^{B-},y_t,Y_t),
\end{align} where $x_t^{G+} = (x^{G+}_{tg})_{g\in\mathcal{G}}$, $x_t^{B+}=(x_{tb}^{B+})_{b\in\mathcal{B}}$, $x_t^{B-}=(x_{tb}^{B-})_{b\in\mathcal{B}}$, $y_t=(y_{tn})_{n\in\mathcal{N}}$ such that $y_t \in \mathcal{Y}$, and $Y_t=(Y_{t,ij})_{(i,j)\in\mathcal{N}\times\mathcal{N}}$. $x_t^{G+}$ and $x_t^{B+}$ represent time $t$ power injections at each node by generators and batteries respectively, while $x_t^{B-}$ represents power sent from the grid to storage devices at each node at time $t$. $y_t$ is a vector describing the power arriving at each node at time $t$ while $Y_{t,ij}$ is the power flow from node $i\in \mathcal{N}$ to $j\in \mathcal{N}$ at time $t$. We have $Y_{t,ij}=-Y_{t,ji}$ if $(i, j) \in \mathcal{E}$, $Y_{t,ij}=0$ if $(i, j) \not\in \mathcal{E}$, and $y_{t,i}=\sum\limits_{j\in \mathcal{N}}Y_{t,ji}$. The decision $x_t$ is also subject to the following constraints for all $t=0,...,T$, $g\in \mathcal{G}$, and $b\in \mathcal{B}$:
\begin{align}
&x_{t,g}^{G+}=\kappa_g^l \text{ if } Z_{0,g}^G=1 \text{ or }  (Z_{t,g}^G=1 \cap Z_{t-1,g}^G=0) \label{D1}\\
&\kappa_g^l Z_{t,g} \leq x_{t,g}^{G+} \leq \kappa_g^u Z_{t,g}, \label{D2}\\
&\kappa_b^l \leq R^B_{t,b} + \eta^-_b x_{t,b}^{B-} - \eta^+_b x_{t,b}^{B+} \leq \kappa_b^u .\label{D4}
\end{align} All vector decisions $x_t$ that satisfy the above constraints form the set of feasible decisions at time $t$, $\mathcal{X}_t(S_t)$.

Constraint \eqref{D1} requires that generators that come online at time $t$ (or start on at time $0$) start at their minimum generation capability. Constraint \eqref{D2} ensures generator power outputs are within minimum and maximum limits. A generator's output must also be $0$ if that generator is scheduled to be off at time $t$ ($Z_{t,g}=0$). Constraint \eqref{D4} disallows storage decisions that violate minimum and maximum capacity limits at each device.

\subsection{Exogenous Information}
\label{Exo Info}


We are interested in the problem of efficiently operating storage devices in the presence of large amounts of offshore wind, which will be modeled as an aggregate quantity. As of this writing (2019), New Jersey has set a goal of 3,500 MW of offshore wind generating capacity.

\begin{figure}
\centering
\includegraphics[width=\columnwidth, height=3.5 in]{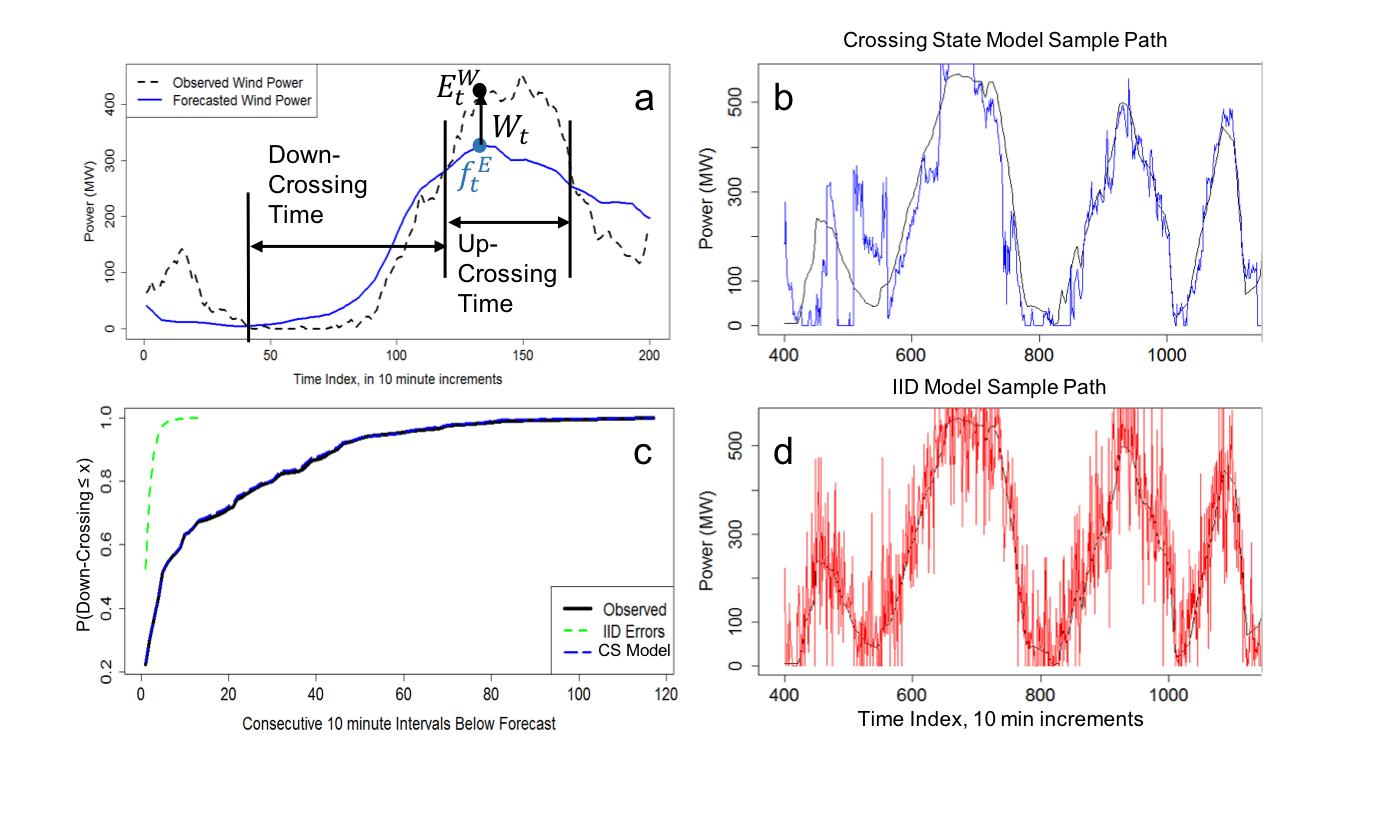}
\caption{a: Examples of crossing times for wind power forecast errors. A single forecast error, $W_t=E^W_t-f_t^E$, is also shown. b and d: Examples of sample paths generated assuming a crossing state model (b) and an IID model (d) for forecast errors. Due to the lack of intertemporal dependencies, the IID model sample paths exhibit far too much variability. c: The empirical cumulative distribution function for crossing times (only down-crossings are shown here) for the crossing state model and IID model over 100 sample paths versus historical distributions for a single wind farm. While the IID model produces crossing times that are far too short compared to reality, this crossing time distribution is very accurately replicated by the crossing state model.}
\label{CrossingTimePic}
\end{figure}

Unit commitment decisions are made leveraging wind power forecasts over the course of the optimization horizon, $\left\lbrace f_t^E \right\rbrace_{t=0}^{T}$. On the real-time, economic dispatch scale, we experience uncertainty in the form of errors from this forecast, $W_t=E_t^W-f_t^E$, which occur on the five minute scale. The forecast error model we use is the univariate crossing state model presented in \cite{durante2017} which accurately replicates the crossing times of wind power forecast errors. Figure \ref{CrossingTimePic}a shows two examples of crossing times. We refer the reader to Appendix A for a detailed description and thorough discussion of the model, but we highlight some important features of the model here:
\begin{hitemize}
\item[1.] Crossing time distributions are explicitly modeled by forming a unique crossing time distribution for each of a small number (typically six) of partially observable states called crossing states. The crossing state is a combination of whether we are above or below the forecast (observable at time $t$) and an aggregated crossing time length variable (i.e. short, medium, long) giving the distribution of time for which wind will remain above or below its forecast. This is hidden to the system operator at time $t$ as we cannot be sure which crossing time bin the system was in until after the wind forecast error switches signs.
\item[2.] For each crossing state, time $t$ forecast errors are aggregated into a small number of discrete states as well (typically three to five), for which distributions of $W_{t+1}$ are formed. The combination of the crossing state and the time $t$ error bin form a partially hidden information state $I_t^W$.
\end{hitemize}

To continue the discussion, we must first define some additional notation for stochastic processes. Let a sample path $\omega$ be given by a sequence of the random variables $W_t$ for $t=1,2, ...,T$: $\omega = (W_1,W_2, ..., W_T)$, and the set of all possible sample path outcomes be given by the sample space $\Omega$. Given a probability space $(\Omega, \mathcal{F},\mathbb{P})$ with sample space $\Omega$, sigma algebra $\mathcal{F}$, probability measure $\mathbb{P}$, and filtration $\sigma\left\lbrace \emptyset, \Omega \right\rbrace =  \mathcal{F}_0 \subset \mathcal{F}_1 \subset \mathcal{F}_2 \subset ... \subset \mathcal{F}_T = \mathcal{F}$, the post-decision information state at time $t$, $I_t^x$, is determined by an $\mathcal{F}_t$-measurable mapping $I_t^x(\omega)$ of sample path $\omega \in \Omega$ to some $I_t^x \in \mathcal{I}_t^x(\Omega)$, where $\mathcal{I}_t^x(\Omega)$ is the set of all possible time $t$ post-decision information states. The post-decision information state contains only the information necessary to determine the distribution $\mathbb{P}(W_{t+1}=w|I_t^x)$ for all $w \in \Omega_{t+1}$, where $\Omega_{t}$ is the outcome space at time $t$. Note that the size of $|\mathcal{I}_t^x(\Omega)|$ is largely dependent on model choice and can explode fairly quickly as more information is needed to transition to the next time period. For example, an autoregressive time series model of order $M$ would have $|\mathcal{I}_t^x(\Omega)|=\prod\limits_{t'=t-M+1}^{t} |\Omega_{t'}|$. Using models that can take on many post-decision information states restricts our ability to index VFAs on these states. This brings us to the next highlight of the crossing state model:
\begin{hitemize}
\item[3.] Given $I_t^W$, we can determine the distribution of $W_{t+1}$. Thus, $E_t^W$ does not need to be included in the post-decision information state as it is not needed to transition to time $t+1$ after a decision has been made. Therefore, we have $I_t^x=I_t^W$ and the set of all possible values for $I_t^x$, the post-decision information state space $\mathcal{I}_t^x(\Omega)$, is the same for all $t$ and relatively compact. It usually consists of somewhere between twelve to thirty discrete states depending on the number of crossing states and error bins that are chosen.
\item[4.] As the information states in the crossing state model are partially hidden, a sample path of wind up to time $t$ may correspond to different information states. Thus, we must track our time $t$ beliefs about the probability that the process is in each information state to aid in decision making at time $t$. This set of beliefs $\mathbb{P}(I_t^x=i)$ for $i \in \mathcal{I}_t^x(\Omega)$ comprises the knowledge state at time $t$, $K_t$. Our belief about the distribution of $W_{t+1}$ given $K_t$ is
\begin{align}
\mathbb{P}(W_{t+1}=w|K_t)=\sum\limits_{i \in \mathcal{I}_t^x(\Omega)}\mathbb{P}(W_{t+1}=w|I_t^x=i)\mathbb{P}(I_t^x=i).
\end{align}
\item[5.] Sample paths produced with this model are realistic and exhibit superior crossing time behavior, especially when compared to a model that assumes errors are independent over time. This is shown in Figure \ref{CrossingTimePic}c and can be subjectively observed by comparing Figure \ref{CrossingTimePic}b and Figure \ref{CrossingTimePic}d.
\end{hitemize}



We also consider an independently and identically distributed (IID) wind power error model. This model exhibits stagewise independence and is thus well suited for use in combination with classic SDDP. This is trained on the same sample path data as the crossing state model; however, VFAs are computed with an algorithm that assumes $W_t$ is distributed according to the same distribution at each time step -- the observed forecast error distribution gathered from data. Under this assumption, no information is needed in the post-decision information state as the distribution of wind at each time step is independent of the history of the exogenous process. Thus, there exists only one possible (trivial) information state and $|\mathcal{I}_t^x(\Omega)|=1$ for all $t$. Note this implies that no information is needed in the knowledge state either as we have probability $1$ of being in the single information state. We compare the effectiveness and robustness of the solution resulting from each wind power model assumption in Section \ref{Numerical Results}.

\subsection{The Transition Function}
\label{Transition Function}

The system transition function, $S_{t+1}=S^M(S_t,x_t,W_{t+1})$, is broken into two steps: the pre-decision to post-decision state transition function given a decision $x_t$, $S_t^x=S^{M,x}(S_t,x_t)$, and the transition function from post-decision state to the next pre-decision state given the arrival of exogenous information $W_{t+1}$, $S_{t+1}=S^{M,W}(S_t^x,W_{t+1})$. To define $S_t^x=S^{M,x}(S_t,x_t)$, we first define $p_t \in \mathbb{R}^{|\mathcal{N}|}$ to be the total power generated at each node. Letting $\mathcal{G}(i)$, $\mathcal{B}(i)$, and $\mathcal{H}(i)$ be the generators, storage devices, and wind farms that map to node $i$, we can calculate the $i$-th component of the vector at time $t$ as
\begin{align}
p_{t,i}=\sum\limits_{g \in \mathcal{G}(i)} x_{t,g}^{G,+}+\sum\limits_{b \in \mathcal{B}(i)} (x_{t,b}^{B,+}-x_{t,b}^{B,-})+\sum\limits_{h \in \mathcal{H}(i)} E_{t,h}^W.
\end{align} $S_t^x=S^{M,x}(S_t,x_t)$ is then given by
\begin{align}
R_{t,b}^{B,x}&=R_{t,b}^{B} + \eta^-_b x_{t,b}^{B-} - \eta^+_b x_{t,b}^{B+} \text{ for } b \in \mathcal{B}, \label{Trans 1}\\
R_{t}^{S,x}&=R_{t}^{S}+\sum\limits_{i \in \mathcal{N}} \max\left\lbrace 0, d_{t,i}-(p_{t,i}+y_{t,i}) \right\rbrace,  \label{Trans 5}\\
I_t^x&= S^{I,x}(I_t) = I_t^W, \label{Trans 2}
\end{align} where $d_{t,i}$ is the demand at node $i$ and $y_{t,i}$ is the power inflow at node $i$. Equation \eqref{Trans 1} alters the post-decision resource state at each storage device according to the charging/discharging decisions made, while Equation \eqref{Trans 5} updates the cumulative shortage state by adding any shortages incurred at time $t$. Equation \eqref{Trans 2} defines the transition between the pre- and post-decision information state. It is a simple function, but explaining why this is the case requires a more in-depth discussion of the crossing state model. These details are left for Appendix A. $S_{t+1}=S^{M,W}(S_t^x,W_{t+1})$ is then defined by
\begin{align}
R_{t+1,b}^{B}&=R_{t,b}^{B,x} \text{ for } b \in \mathcal{B}, \label{Trans 3}\\
R_{t+1}^{S}&=R_{t}^{S,x}, \label{Trans 6}\\
(K_{t+1},I_{t+1})&=S^{K,W}(K_t,I_t^x,W_{t+1}), \label{Trans 4}
\end{align} where we assume the battery state and shortage state remain unchanged from post- to pre- decision state (Equation \eqref{Trans 3} and \eqref{Trans 6}).  Following an observation $W_{t+1}=w$, we update our knowledge and information states using the Bayesian updating function $(K_{t+1},I_{t+1})=S^{K,W}(K_t,I_t^x,W_{t+1})$. We direct the reader to Appendix A for specifics on how this is formulated using the crossing state model. Note, though, that under the IID wind model assumption, $(K_{t+1},I_{t+1})=(K_t,I_t^x)=(K_t, I_t)$ for all $t$ as there is a single information state (pre- and post-decision), and $\mathbb{P}(W_{t+1}=w|K_t)$ can be simplified to $\mathbb{P}(W_{t+1}=w)$.

\subsection{The Objective Function}
\label{Objective Function}

At time $t$, the cost of generation and utilizing storage is
\begin{align}
C^G_t(S_t,x_t)=\sum\limits_{g \in \mathcal{G}} c_{t,g}^G x_{t,g}^{G,+} +\sum\limits_{b \in \mathcal{B}} c_{t,b}^B (x_{t,b}^{B,+}-x_{t,b}^{B,-}).
\end{align}

We then introduce penalties to discourage decisions that result in shortages, excess in generation, and violation of line constraints at each time period. These are given respectively by
\begin{align}
C^S_t(S_t,x_t)&=\theta^S \sum\limits_{i \in \mathcal{N}} \max\left\lbrace 0, d_{t,i}-(p_{t,i}+y_{t,i}) \right\rbrace,\\
C^E_t(S_t,x_t)&=\theta^E \sum\limits_{i \in \mathcal{N}} \max\left\lbrace 0, (p_{t,i}+y_{t,i})-d_{t,i} \right\rbrace,\\
C^L_t(S_t,x_t)&=\theta^L \sum\limits_{(i, j) \in \mathcal{E}} \max\left\lbrace 0, |Y_{t,ij}|-\ell_{ij} \right\rbrace,
\end{align} where $\theta^S$, $\theta^E$, and $\theta^L$ are nonnegative parameters set by the system operator.

The shortage state threshold penalty is applied at time $T$ as follows,
\begin{align}
C^P_t(S_t,x_t) = \begin{cases}
\theta^P \max\left\lbrace 0, R_t^S-\theta^C \right\rbrace \qquad &\text{ if } t=T, \\
0 \qquad &\text{ otherwise}. 
\end{cases}
\end{align}The nonnegative parameter $\theta^C$ can be viewed as the ``maximum acceptable cumulative shortages," determined by the system operator, while nonnegative parameter $\theta^P$ controls how much shortages in excess of this are penalized. The parameter $\theta^P$ can be tuned by a system operator to manage the risk-reward trade-off as a higher penalty will likely lead to more robust, but more expensive, solutions. It is important to note that this threshold penalty is convex with respect to $R_T^S$. As SDDP requires convexity in the resource state, this is an appropriate choice of penalty for the shortage state, as opposed to, for example, a step function penalty.

The total cost is then given by the sum of the generation costs and penalty functions,
\begin{align}
C_t(S_t,x_t)=C^G_t(S_t,x_t)+C^S_t(S_t,x_t)+C^E_t(S_t,x_t)+C^L_t(S_t,x_t)+C^P_t(S_t,x_t).
\end{align} 

We aim to operate the system at minimize cost over a finite-horizon in five-minute intervals over the course of one day. Thus our stochastic optimization problem has the objective function
\begin{equation}
\label{Obj}
\operatornamewithlimits{min}\limits_{\pi \in \Pi}\mathbb{E} \left[\sum\limits_{t=0}^T C(S_t,X^{\pi}_t(S_t))|S_0\right],
\end{equation} where $T=288$ and the system transition functions are given by $S_t^x=S^{M,x}(S_t,x_t)$ and $S_{t+1}=S^{M,W}(S_t^x,W_{t+1})$. $X^{\pi}_t(S_t)$ maps the state $S_t$ to an action using policy $\pi$, and we are minimizing over the set of all admissible policies $\pi \in \Pi$.

\section{Development of the SDDP Algorithm}
\label{SDDP Algo}

Stochastic dual dynamic programming is an algorithm originating from the stochastic programming community designed to solve multistage stochastic linear programs of the form 
\begin{align}
\label{StochProg}
\min\limits_{\substack{A_0 x_0=b_0 \\ x_0\geq 0}}\left\langle c_0, x_0 \right\rangle+ \mathbb{E}_1 \left[\min\limits_{\substack{B_0 x_0+A_1 x_1=b_1 \\ x_1\geq 0}}\left\langle c_1, x_1 \right\rangle + \mathbb{E}_2 \left[ \cdots +\mathbb{E}_T \left[\min\limits_{\substack{B_{T-1} x_{T-1}+A_T x_T=b_T \\ x_T\geq 0}}\left\langle c_T, x_T \right\rangle \right] \cdots \right] \right],
\end{align} where $W_t=(A_t,B_t, b_t, c_t)$ is the random information at time $t$. $A_t$ and $B_t$ are $\mathcal{F}_t$-measurable random matrices and $b_t$ and $c_t$ are $\mathcal{F}_t$-measurable random vectors. The initial state $S_0=(A_0,B_0, b_0, c_0)$ is an exception as it is assumed to be deterministic.

We now show how these multistage stochastic linear programs can be written using dynamic programming notation. Recall the dynamic pre- and post-decision state variables are $S_t=(R_t,I_t,K_t)$ and $S_t^x=(R_t^x,I_t^x,K_t)$ respectively where the information and knowledge states contain all the information necessary to model and control the system from time $t$ onward that is not contained in the resource states. By defining $R_t=B_{t-1}x_{t-1}-b_t$, $R_t^x=B_t^x x_t$, $C(S_t,x_t)=\left\langle c_t, x_t \right\rangle$, and the set of feasible decisions
\begin{align}
\mathcal{X}_t(S_t)=\begin{cases}
x_t \in \mathbb{R}^{\dim(c_t)}: A_t x_t = b_t \qquad &\text{ if } t=0, \\
 x_t \in \mathbb{R}^{\dim(c_t)}: B_{t-1}x_{t-1} + A_t x_t = b_t \qquad &\text{ if } t>0, 
\end{cases}
\end{align} we can rewrite the optimization problem in \eqref{StochProg} as
\begin{align}
\label{DynProg}
\min\limits_{x_0 \in \mathcal{X}_0(S_0)} C(S_0,x_0)+ \mathbb{E}_1 \left[\min\limits_{x_1 \in \mathcal{X}_1(S_1)} C(S_1,x_1) + \mathbb{E}_2 \left[ \cdots +\mathbb{E}_T \left[\min\limits_{x_T \in \mathcal{X}_T(S_T)} C(S_T,x_T) \right] \cdots \right] \right].
\end{align} As the problem is stochastic, the optimal solution is not a vector of decisions, but instead a policy $\pi$ that maps a state $S_t$ to an action $x_t=X^{\pi}_t(S_t)$. We thus aim to minimize over the set of admissible policies, which is the set of $\mathcal{F}_t$-measurable functions that map states to actions, and the optimization problem in \eqref{DynProg} can be written in the form in \eqref{Obj}.

The optimal policy may be found by finding exact value functions, or expected cost-to-go functions, for every possible system state. These are given by Bellman's equation for finite horizon problems, \begin{align}
\label{ExactVal}
V_t^*(S_t)=\min\limits_{x_t \in \mathcal{X}_t(S_t)} \left(C(S_t,x_t)+\mathbb{E}\left[V^*_{t+1}(S_{t+1})|S_t,x_t\right]\right),
\end{align}
and, once these are found, the optimal policy is given by \begin{align}
\label{ExactPol}
X_t^{*}(S_t)=\argmin\limits_{x_t \in \mathcal{X}_t(S_t)} \left(C(S_t,x_t)+\mathbb{E}\left[V^*_{t+1}(S_{t+1})|S_t,x_t\right]\right).
\end{align} When finding exact value functions for all possible system states is computationally intractable due to the curses of dimensionality (as is the case in our problem), we can instead rely on approximations of these value functions and use a VFA-based policy, \begin{align}
\label{ApproxPol}
X_t^{\pi}(S_t)=\argmin\limits_{x_t \in \mathcal{X}_t(S_t)} \left(C(S_t,x_t)+\mathbb{E}\left[\bar{V}_{t+1}(S_{t+1})|S_t,x_t\right]\right),
\end{align} where $\bar{V}_{t+1}(S_{t+1})$ is some approximation of the value of the downstream states.

Alternatively, we can fit value functions instead to the post-decision state variable $S_t^x$. The resulting post-decision state-based VFA policy is given by \begin{align}
\label{ApproxPolPDS}
X_t^{\pi}(S_t)=\argmin\limits_{x_t \in \mathcal{X}_t(S_t)} \left(C(S_t,x_t)+\bar{V}^x_{t}(S^x_{t})\right),
\end{align}
where $\bar{V}^x_{t}(S^x_{t})$ serves as an approximation of $V^{x,*}_{t}(S^x_{t})=\mathbb{E}\left[V_{t+1}^*(S_{t+1})|S_t^x\right]$. Noting that $S_t^x=(R_t^x,I_t^x,K_t)$ in this problem, where $K_t$ contains time $t$ beliefs giving  $\mathbb{P}(I_t^x=i)$ for $i \in \mathcal{I}_t^x(\Omega)$, this approximation can be written $\bar{V}^x_{t}(S^x_{t})=\sum\limits_{i \in \mathcal{I}_t^x(\Omega)}\bar{V}^x_t(R_t^x,I_t^x=i)\mathbb{P}(I_t^x=i)$ and we can rewrite Equation \eqref{ApproxPolPDS} as
\begin{align}
\label{ApproxPolPDS2}
X_t^{\pi}(S_t)=\argmin\limits_{x_t \in \mathcal{X}_t(S_t)} \left(C(S_t,x_t)+\sum\limits_{i \in \mathcal{I}_t^x(\Omega)} \bar{V}^x_t(R_t^x,I_t^x=i)\mathbb{P}(I_t^x=i)\right).
\end{align} The challenge is now fitting these value functions, $\bar{V}^x_{t}(R_t^x,I_t^x)$, in a setting with a high-dimensional resource state and a hidden Markov information state. The remainder of this section describes the development of the variant of SDDP used to accomplish this. Section \ref{Classic SDDP} reviews classical SDDP methods, Section \ref{StagewiseDependentSDDP} details how to incorporate a hidden Markov information state, and Section \ref{IS Section} describes an importance sampling technique used to speed up convergence without sacrificing solution quality.

\subsection{Classic SDDP Algorithms}
\label{Classic SDDP}

Classic SDDP, introduced by \cite{pereira1991multi}, is an iterative two pass procedure that alternates between forward and backward passes to fit post-decision state VFAs. The classic algorithm is limited in that it requires interstage independence for all random quantities, implying that the cardinality of the post-decision information state space at each time $t$, $|\mathcal{I}_t^x(\Omega)|$, is one. Under this assumption, we simplify our policy to $X_t^{\pi}(S_t)=\argmin\limits_{x_t \in \mathcal{X}_t(S_t)} \left(C(S_t,x_t)+\bar{V}^x_{t}(R_t^x) \right)$ and fit value functions to only the post-decision resource state. SDDP overcomes the high-dimensionality of the resource state by using piecewise linear value function approximations, $\bar{V}^x_{t}(R^x_{t})$, formed as the maximum of a set of cutting hyperplanes known as Benders cuts, which are lower-bounding convex outer approximations to a convex exact value function $V_t^{x,*}(R_t^x)$. Convexity of the value function in the resource state (a property of our energy storage problem) is assumed as it is necessary for SDDP. A high-level outline of the classic SDDP algorithm is given in Algorithm \ref{alg:HighLevel}. This is accompanied by notation and details that appear in the subsequent paragraph.

During a forward pass of the algorithm at iteration $k$, the VFAs from the previous iteration, $\bar{V}_{t}^{x,k-1}(R_t^x)$, form the policy which determines the resource points $R_t^{x,k}$ at each time step $t$. These lower-bounding VFAs have the form 
\begin{align}
\label{benders}
\bar{V}_t^{x,k-1}(R_t^x)=\max\limits_{i \leq k-1}\left\lbrace \alpha_t^i +\left\langle \beta_t^i, R_t^x-R_t^{x,i} \right\rangle \right\rbrace.
\end{align} The decision at each time step in the forward pass is determined by solving the linear program
\begin{align}
\label{LP}
x_{t}^k \in \argmin\limits_{x_{t} \in \mathcal{X}_{t}(S_t)} \left\lbrace C(S_{t},x_{t}) + \bar{V}_{t}^{x,k-1}(R_{t}^x) \right\rbrace,
\end{align}
and the subsequent post-decision resource point is determined by altering the pre-decision state resource vector according to the decision $x_t^k$: $R_t^{x,k}=S^{M,x}(S_t,x_t)$.

Subsequently, in the backward pass of the same iteration, these resource points are utilized to create new Benders cuts that are lower bounds to $V_t^{x,*}(R_t^x)$ and update VFAs in the process. The Benders cut, $h^k_t(R_t^x)$, at time $t$ and iteration $k$ has the form
\begin{align}
h^k_t(R_t^x)= \ubar{V}_{t+1}^{x,k}(R_t^{x,k})+\left\langle\beta_t^k,R_t^x-R_t^{x,k}\right\rangle,
\end{align} where $ \ubar{V}_{t+1}^{x,k}(R_t^{x,k})$ is a lower bound for $V_{t}^{x,*}(R_t^{x,k})$. We now must find $\ubar{V}_{t+1}^{x,k}(R_t^{x,k})$ and $\beta_t^k$ to construct the cut. Assuming a discrete outcome space $\Omega_{t+1}$ and defining
\begin{align}
\ubar{V}_{t+1}^{x,k} (R_{t}^x,w)=\min\limits_{x_{t+1} \in \mathcal{X}_{t+1}(S_{t+1}(w))} \left\lbrace C(S_{t+1}(w),x_{t+1}) + \bar{V}_{t+1}^{x,k}(R_{t+1}^x) \right\rbrace,
\end{align} for $w \in \Omega_{t+1}$, where $\bar{V}_{t+1}^{x,k}(R_{t+1}^x)$ is a lower bound for $V_{t+1}^{x,*}(R_{t+1}^x)$ of the form given by Equation \eqref{benders}, we can find a lower bound for the value at any resource point $R_{t}^x$ by taking the expectation
\begin{align}
\label{expVal}
\ubar{V}^{x,k}_{t+1}(R_t^x)=\sum\limits_{w \in \Omega_{t+1}} \mathbb{P}(w)\ubar{V}_{t+1}^{x,k} (R_t^{x},w).
\end{align} Using $R_{t}^{x,k}$ as the argument, we can find $\ubar{V}_{t+1}^{x,k}(R_{t}^{x,k})$, which also serves as $\alpha_t^k$ in Equation \eqref{benders}. Furthermore, letting $\beta_t^k(w)$ be a subgradient of $\ubar{V}_{t+1}^{x,k} (R_{t}^{x,k},w)$ with respect to $R_t^x$ for each $w \in \Omega_{t+1}$, we can choose a subgradient of Equation \eqref{expVal} at $R_t^{x,k}$ and set
\begin{align}
\label{subGrad}
\beta_t^k \in \delta_R \ubar{V}^{x,k}_{x,t+1}(R_t^{x,k}) = \sum\limits_{w \in \Omega_{t+1}} \mathbb{P}(w)\beta_t^k(w).
\end{align} This completes the cut $h_t^k(R_t^x)$. Finally, to update the VFA at time $t$, we take the maximum of the new cut and the VFA from the previous iteration,
\begin{align}
\bar{V}_t^{x,k}(R_t^x)=\max\left\lbrace \bar{V}_t^{x,k-1}(R_t^x),h_t^k(R_t^x) \right\rbrace,
\end{align} which is of the form given by Equation \eqref{benders}.

\begin{algorithm}
\footnotesize
\caption{A Basic Outline of Classic SDDP}
\begin{algorithmic}[1]
\STATE Initialize system state $S_0$, which includes $R_0^{k}=R_{-1}^{x,k}$ for each iteration $k$, as well as value function approximations $\bar{V}_t^{x,0}(R_t^x)= -\infty$ for $t<T$, and $\bar{V}_T^{x,k}(R_T^x)= V_T^{x,*}(R_T^x)$ for all $k$.
\FOR{iteration $k=0,...,K$}
\STATE \textit{Forward Pass, used to determine resource points visited in iteration $k$:}
\STATE Generate scenario $\omega \in \Omega$.
\FOR {$t = 0, ..., T$}
\STATE Based on the system state and current VFAs, solve $x_t ^k=\argmin\limits_{x_t \in \mathcal{X}_t(S_{t}(\omega))}\left\lbrace C(S_t(\omega),x_t)+ \bar{V}_t^{x,k-1}(R_{t}^{x})\right\rbrace$.
\STATE $R_t^{x,k}$ and then $S_{t+1}(\omega)$ are determined by the decision $x_t ^k$ and the system transition functions.
\ENDFOR
\STATE \textit{Backward Pass, used to update VFAs in iteration $k$:}
\FOR {$t = T, ..., 1$}
\FOR {each discrete outcome $w \in \Omega_t$}
\STATE Let $\ubar{V}_{t}^{x,k} (R_{t-1}^{x},w) = \min\limits_{x_t \in \mathcal{X}_t(S_t(w))} \left\lbrace C(S_t(w),x_t) +  \bar{V}_t^{x,k}(R_t^{x}) \right\rbrace$.
\STATE At $R_{t-1}^{x,k}$, find $\ubar{V}_t^{x,k} (R_{t-1}^{x,k},w)$ and subgradient $\beta_{t-1}^k(w) \in \delta_R \ubar{V}^{x,k}_{t-1}(R_{t-1}^{x,k})$.
\ENDFOR
\STATE Find the hyperplane parameters $\ubar{V}^{x,k}_{x,t}(R_{t-1}^{x,k})$ and $\beta_{t-1}^k$ at $R_{t-1}^{x,k}$ by taking the expectation over  $\Omega_{t+1}$:
\STATE $\ubar{V}^{x,k}_t(R_{t-1}^{x,k})=\sum\limits_{w \in \Omega_{t+1}} \mathbb{P}(w)\ubar{V}_{t}^{x,k} (R_{t-1}^{x,k},w)$, $\beta_{t-1}^k=\sum\limits_{w \in \Omega_{t+1}} \mathbb{P}(w)\beta_{t-1}^k(w)$.
\STATE These define the time $t-1$ Benders cut in iteration $k$, $h_{t-1}^k(R_{t-1}^x)=\ubar{V}^{x,k}_{x,t}(R_{t-1}^{x,k})+\left\langle \beta_{t-1}^k, R_{t-1}^x-R_{t-1}^{x,k}\right\rangle$.  
\STATE The VFA at time $t-1$ is now the maximum of a set of $k$ Benders cuts. This can be written in a recursive updating form as follows: $\bar{V}_{t-1}^{x,k}(R_{t-1}^x)=\max\left\lbrace \bar{V}_{t-1}^{x,k-1}(R_{t-1}^x), h_{t-1}^k(R_{t-1}^x) \right \rbrace$.
\ENDFOR
\ENDFOR
\end{algorithmic}
\label{alg:HighLevel}
\end{algorithm}

These value functions $\bar{V}_t^{x,k}(R_t^x)$ are lower bounding approximations of $V_t^{x,*}(R_t^x)$ and are monotonically increasing with respect to iteration $k$. They are shown in \cite{philpott2008convergence} to converge to the optimal value $V_t^{x,*}(R_t^x)$ in a finite number of iterations with probability 1 under the following assumptions:

\begin{itemize}
\item[A1.] Random quantities appear only on the right-hand side of the linear constraints in each stage.
\item[A2.] The set $\Omega_t$ of random outcomes in each stage $t = 2, 3,..., T$ is discrete and finite.
\item[A3.] Random quantities in different stages are independent.
\item[A4.] The feasible region of the linear program in each stage is non-empty and bounded.
\end{itemize}

As written, calculating the subgradients necessary to form a new Benders cut at each time step and iteration (Equations \eqref{LP}-\eqref{subGrad}) entails finding dual solutions to the linear program defined by equation \eqref{LP} for each $w \in \Omega_{t+1}$. As this is often the computational bottleneck in the algorithm, previous works have attempted to accelerate this step. For example, cut selection can shorten the time required to solve each linear program by reducing the number of cuts involved in the linear programs at each stage \citep[see][]{de2015improving}. Alternatively, one may seek to reduce the number of linear programs solved at each stage. \cite{higle1991stochastic}, \cite{chen1999convergent}, and \cite{hindsberger2014resa} all present algorithms which form cuts based on a sampled subset of the outcome space $\tilde{\Omega}_{t+1} \subseteq \Omega_{t+1}$. An example of a sample-based cut calcualtion algorithm related to our work appears in \cite{philpott2008convergence}. In the backward pass, $\tilde{\Omega}_{t+1}$ is determined by sampling outcomes from some distribution $\mathbb{Q}^k_t$. Next, extreme-point dual solutions for $w \in \tilde{\Omega}_{t+1}$ are found and stored at each iteration $k$. Then, for $w \in \Omega_{t+1}$, one finds the best dual solution of those stored to form an approximate lower bounding hyperplane at $w$ parameterized by $\hat{\ubar{V}}_{t+1}^{x,k}(R_t^{x,k},w)$ and $\hat{\beta}_t^k(w)$. The cut is then calculated by taking an expectation over the entire sample space in a manner similar to equations \eqref{expVal} and \eqref{subGrad}, but utilizing $\hat{\ubar{V}}_{t+1}^{x,k}(R_t^{x,k},w)$ and $\hat{\beta}_t^k(w)$ in place of $\ubar{V}_{t+1}^{x,k}(R_t^{x,k},w)$ and $\beta_t^k(w)$ respectively.


\subsection{SDDP with Hidden Markov Models, Quadratic Regularization, and Sampling in the Backward Pass}
\label{StagewiseDependentSDDP}

Note that one of the assumptions made by classical SDDP is that of stagewise independence for all random quantities, greatly limiting our choice of stochastic models. IID models are often chosen for the stochastics as they are a straightforward choice that satify this requirement. Under certain circumstances (i.e. interstage dependence only occurs in the right hand side constraint vectors) we can extend the formulation to incorporate autoregressive moving average (ARMA) models. This extension requires additional state variables to track the history of the process and an increase in the dimensionality of the value function approximation that depends on the order of the model. Besides the downside of potentially much slower convergence (though cut sharing algorithms such as those from \cite{infanger1996cut} help circumvent this issue), often a stochastic process cannot be properly modeled by a standard ARMA process. This is the case, for example, with solar power which may be best described by regime switching models \citep[as in][]{shakya2016solar}. In our case we aim to model wind power with a hidden Markov model, and thus need an algorithm that will allow for general forms of interstage dependence. Various methods accomplish this, including the approximate dual dynamic programming algorithms seen in \cite{lohndorf2013optimizing} and \cite{lohndorf2015optimal}.

Most applicable to our problem, however, is the algorithm from \cite{asamov2015regularized}, which accommodates general Markov uncertainty (beyond ARMA processes) in the stochastics by fitting a set of Benders cuts to each possible post-decsision information state at each time step. This is effectively a lookup table representation of the value function, since we have a different set of cuts for each possible discretized value of the information state. This work assumes that the post-decision state is fully observable. Assuming $|\mathcal{I}_t^x(\Omega)|$ is finite, VFAs are indexed on the post-decision information state $I_t^x \in \mathcal{I}_{t}^x(\Omega)$, and are of the form
\begin{align}
\label{benders2}
\bar{V}_t^{x,k}(R_t^x,I_t^x)=\max\limits_{i \leq k}\left\lbrace \alpha_t^i(I_t^x) +\left\langle \beta_t^i(I_t^x), R_t^x-R_t^{x,i} \right\rangle \right\rbrace.
\end{align}
In the forward pass, $W_t$ is sampled from the Markov process at each time step, and the post-decision resource points visited in iteration $k$, $R_t^{x,k}$, are then determined by following the policy given by equation \eqref{ApproxPolPDS} using the VFAs from iteration $k-1$, $\bar{V}_t^{x,k-1}(R_t^x,I_t^x)$. The value functions for each $I_t^x \in \mathcal{I}_{t}^x(\Omega)$ are then updated in the backward pass by taking the maximum of $\bar{V}_t^{x,k-1}(R_t^x,I_t^x)$ and the Benders cut $h_t^k(R_t^x,I_t^x)= \alpha_{t}^k(I_{t}^x)+\left\langle\beta_t^k(I_t^x),R_t^x-R_t^{x,k}\right\rangle$ constructed in iteration $k$. Conditional probabilities can now be used to form the intercept and slope vectors in the following manner,
\begin{align}
\alpha_t^k(I_t^x)\leftarrow\ubar{V}_{t+1}^{x,k} (R_t^{x},I_t^x)=\sum\limits_{w \in \Omega_{t+1}} \ubar{V}_{t+1}^{x,k} (R_t^{x},w) \mathbb{P}(w|I_t^x),
\end{align} and
\begin{align}
\beta_t^k(I_t^x) \in \delta_R \ubar{V}^{x,k}_{t+1}(R_t^{x,k},I_t^x).
\end{align}

We extend the algorithm from \cite{asamov2015regularized} to accommodate hidden Markov stochastic models. As the information states are now partially unobservable, the post-decision resource points visited in iteration $k$, $R_t^{x,k}$, are determined in the forward pass by following the policy given by equation \eqref{ApproxPolPDS2}. We also construct our sample path using the hidden Markov model in the forward pass. The backward pass is then altered slightly by defining, for $w \in \Omega_t$,
\begin{align}
\ubar{V}_t^{x,k} (R_{t-1}^x,w)=\min\limits_{x_t \in \mathcal{X}_t(R_{t-1}^x,I_t(w),K_t(w))} \left\lbrace C(S_t(w),x_t) + \sum\limits_{i \in \mathcal{I}_t^x(\Omega)}\bar{V}_t^{x,k}(R_t^x,I_t^x=i) P(I_t^x=i|w)) \right\rbrace,
\end{align} to account for the fact that an observation of an element of $\Omega_t$ may correspond to different information states. \cite{dowsonapartially} also developed an algorithm for partially observable multiperiod stochastic programs. The implementation differs in that it focuses on finding the value of an initial resource state and underlying state distribution. This is reflected in the algorithm as the forward pass concludes upon a sampled change of the underlying state, upon which a backward pass is performed. Only information states that could possibly be reached by the random realization in the forward pass are visited in the backward pass. Conversely, this paper considers time-dependency in a finite horizon problem and, on each iteration, performs one full forward pass to time $T$. The indexing of information states by time $t$ allows for the fact that the underlying state distribution and transition probabilities may be time-dependent. We also visit all possible information states at each time $t$ in the backward pass, utilizing the resource points visited in the forward pass. 

Another contribution of \cite{asamov2015regularized} was the development of a quadratic regularization technique to accelerate convergence that remains computationally tractable in problems with long time horizons. A regularization term $\frac{\rho^k}{2}\left\langle R_t^x - \bar{R}_t^{x,k-1}, Q_t(R_t^x - \bar{R}_t^{x,k-1})\right\rangle$ is added to equation \eqref{LP} when determining the policy in the forward pass (see Line \ref{regLine} of Algorithm \ref{alg:SDDP}) to penalize large deviations from the incumbent solution (the solution in the previous iteration). $\left\lbrace \rho^k \right\rbrace$ is a regularization sequence such that $\rho^k \geq 0$ for all $k$ and $\lim\limits_{k \rightarrow \infty} \rho^k = 0$. $Q_t$ is a positive semi-definite matrix (typically diagonal) which may be used to address differences in scale across different components of $R_t^x$ by weighting the deviation in each dimension differently. This effectively guides the solution to regions of the value function domain where more information is known and numerical convergence is reached in fewer iterations. Related regularization methods have been used in SDDP previously in \cite{ruszczynski1993regularized}, \cite{higle1994finite}, \cite{morton1996enhanced}, \cite{ruszczynski1997accelerating}, and \cite{sen2014multistage}, for example. However, these methods utilize scenario trees (which experience exponential growth) and are not well suited for the application considered in this paper as the time horizon $T$ is relatively large. One subtlety about the use of regularization in this problem is worth pointing out as it is important to solution quality. Applying regularization on the shortage state will discourage potential reductions in shortages between iterations. As this is undesirable, $Q_t$ is designed such that no penalty is placed on changes in the $R^{S,x}_t$ dimension. 


To further accelerate convergence, we apply sampling in the backward pass, as discussed in Section \ref{Classic SDDP}, to reduce the number of linear programs that must be solved at each stage. Algorithm \ref{alg:SDDP} outlines our solution algorithm, where the sampling method is left unspecified for now. As a placeholder, the sampling distribution at each time step is denoted by a general distribution $\mathbb{Q}^k_t$. Based on successive samples and observations of value functions at each iteration, we may choose to update our sampling distribution $\mathbb{Q}^k_t$. We will discuss alternatives for forming and updating $\mathbb{Q}^k_t$ to sample the outcome space in Section \ref{IS Section}. Additionally, for $w \in \Omega_t$, let $\mathcal{L}(w|I_{t-1}^x)=\mathbb{P} (w|I_{t-1}^x)/\mathbb{Q}^k_t(w)$ be the likelihood ratio giving the ratio of the probability of $w$ occurring under $ \mathbb{P} (w|I_{t-1}^x)$ to the probability of sampling $w$ according to $\mathbb{Q}^k_t$.

\begin{algorithm}
\scriptsize
\caption{SDDP with (Hidden) Markov Models, Quadratic Regularization, and Sampling in the Backward Pass}
\begin{algorithmic}[1]
\STATE Choose $Q_t \succeq 0$, $t=0,...,T$, and define regularization coefficient sequence $\left\lbrace \rho^k \right\rbrace$.
\STATE Let $\bar{V}^{x,k}_T(R_T^x,I_T^x)=V_T^{x,*}(R_T^x,I_T^x)$ for $k=0,...,K$ and $I_T^x\in \mathcal{I}_T^x(\Omega)$.
\STATE Let $\bar{V}^{x,0}_t(R_t^x,I_t^x)=-\infty$ for $I_t^x \in \mathcal{I}_t^x(\Omega)$ and $t=0,...,T-1$.
\STATE Initialize \textit{sampling distributions} $\mathbb{Q}^0_t$ for $t=1,...,T$.
\STATE Initialize $(R^{x,k}_{-1},I_0,K_0) \leftarrow S_0$, $k=0,...,K$.
\FOR{$k=0,...,K$}
\STATE \textit{Forward Pass:}
\STATE Generate sample path $\omega \in \Omega$ using the (hidden) Markov stochastic process.
\FOR {$t = 0, ..., T$}
\IF {$k=0$}
\STATE $x_t ^k=\argmin\limits_{x_t \in \mathcal{X}_t(R_{t-1}^{x,k},I_t(\omega), K_t(\omega))}\left\lbrace C(S_t(\omega),x_t)) \right\rbrace$.
\ELSE
\IF {$t<T$}
\STATE {\begin{eqnarray*}
&x_t ^k=\argmin\limits_{x_t \in \mathcal{X}_t(R_{t-1}^{x,k},I_t(\omega),K_t(\omega))} \Bigg\lbrace C(S_t(\omega),x_t)&+\sum\limits_{i \in \mathcal{I}_t^x(\Omega)} \bar{V}_t^{x,k-1}(R_t^x,I_t^x=i)\mathbb{P}(I_t^x=i|\omega)\\
&&+\frac{\rho^k}{2}\left\langle R_t^x - \bar{R}_t^{x,k-1}, Q_t(R_t^x - \bar{R}_t^{x,k-1})\right\rangle \Bigg\rbrace.
\end{eqnarray*}\label{regLine}}
\ELSE
\STATE \begin{eqnarray*}&x_t ^k=\argmin\limits_{x_t \in \mathcal{X}_t(R_{t-1}^{x,k},I_t(\omega),K_t(\omega))}\left\lbrace C(S_t(\omega),x_t)) + \sum\limits_{i \in \mathcal{I}_t^x(\Omega)} \bar{V}_t^{x,k-1}(R_t^x,I_t^x=i)\mathbb{P}(I_t^x=i|\omega)\right\rbrace.
\end{eqnarray*}
\ENDIF
\ENDIF
\STATE \parbox[t]{\dimexpr\linewidth-\algorithmicindent}{Set $S_t^{x}(\omega)=(R_t^{x,k}, I_t^x(\omega), K_t(\omega)) \leftarrow S^{M,x}(S_t(\omega),x_t^k)$ and $S_{t+1}(\omega)=(R_t^{k}, I_t(\omega),K_t(\omega))\leftarrow S^{M,W}(S_t^{x}(\omega),W_{t+1})$.\strut}
\ENDFOR
\STATE \textit{Backward Pass}
\FOR {$t = T, ..., 1$}
\STATE Sample $\tilde{\Omega}_t$, a small subset of $\Omega_t$, from the current \textit{sampling distribution} $\mathbb{Q}^{k}_t$.
\FOR {$w \in \tilde{\Omega}_t$} \label{LineD}
\STATE Define $\ubar{V}_t^{x,k} (R_{t-1}^x,w)=\min\limits_{x_t \in \mathcal{X}_t(R_{t-1}^x,I_t(w),K_t(w))} \left\lbrace C(S_t(w),x_t) + \sum\limits_{i \in \mathcal{I}_t^x(\Omega)} \bar{V}_t^{x,k}(R_t^x,I_t^x=i)\mathbb{P}(I_t^x=i|w) \right\rbrace$.  \label{LineE}
\STATE Select $\ubar{\beta}_t^k (w) \in \partial_{R_{t-1}^x} \ubar{V}_t^{x,k} (R_{t-1}^{x,k},w)$.  \label{LineF}
\ENDFOR
\FOR{\textbf{all} $I_{t-1}^x \in \mathcal{I}_{t-1}^x(\Omega)$}
\STATE $\alpha_{t-1}^k(I_{t-1}^x)\leftarrow \sum\limits_{w \in \tilde{\Omega}_t} \ubar{V}_t^{x,k} (R_{t-1}^{x,k},w) \mathbb{P}(w|I_{t-1}^x) \mathcal{L}(w|I_{t-1}^x)$; $\beta_{t-1}^k(I_{t-1}^x)\leftarrow \sum\limits_{w \in \tilde{\Omega}_t}\ubar{\beta}_t^k (w) \mathbb{P}(w|I_{t-1}^x)\mathcal{L}(w|I_{t-1}^x)  $. \label{LineA}
\STATE $h_{t-1}^k(R_{t-1}^x, I_{t-1}^x)=\alpha_{t-1}^k(I_{t-1}^x)+\left\langle \beta_{t-1}^k(I_{t-1}^x), R_{t-1}^x-R_{t-1}^{x,k}\right\rangle$.  \label{LineB}
\STATE $\bar{V}_{t-1}^{x,k}(R_{t-1}^x, I_{t-1}^x)=\max\left\lbrace \bar{V}_{t-1}^{x,k-1}(R_{t-1}^x, I_{t-1}^x), h_{t-1}^k(R_{t-1}^x, I_{t-1}^x) \right \rbrace$.  \label{LineC}
\ENDFOR
\STATE{Update the \textit{sampling distribution} $\mathbb{Q}^{k+1}_t$ based on $\mathbb{Q}^{k}_t$ and observed values $\ubar{V}_t^{x,k} (R_{t-1}^{x,k},w)$ for $w \in \tilde{\Omega}_t$.}
\ENDFOR
\STATE $\ubar{V}_0^{x,k} \leftarrow \min\limits_{x_0 \in \mathcal{X}_0(S_0)} \left\lbrace C(S_0,x_0) +\sum\limits_{i \in \mathcal{I}_0^x(\Omega)} \bar{V}_t^{x,k}(R_0^x,I_0^x=i)\mathbb{P}(I_0^x=i) \right\rbrace $.
\STATE $\bar{R}_t^{x,k} \leftarrow R_t^{x,k}$ for $t=0,...,T-1$.
\ENDFOR
\end{algorithmic}
\label{alg:SDDP}
\end{algorithm}

\subsection{Risk-Directed Importance Sampling}
\label{IS Section}

In this section we describe two different methods for determining the sampling distributions $\mathbb{Q}^k_t$. First, however, notice that in order to form non-zero hyperplane parameter vectors $\alpha_{t-1}^k(I_{t-1}^x)$ and $\beta_{t-1}^k(I_{t-1}^x)$ for each $I_{t-1}^x \in \mathcal{I}_{t-1}^x(\Omega)$ in Algorithm \ref{alg:SDDP}, we must ensure that $\mathbb{P}(w|I_{t-1}^x) >0$ for at least one element $w \in \tilde{\Omega}_t$ for every $I_{t-1}^x \in \mathcal{I}_{t-1}^x(\Omega)$.

One simple method to guarantee this is to sample outcomes directly from $\mathbb{P}(W_t=w|I_{t-1}^x)$, the nominal probability distribution given $I_{t-1}^x$, $m$ times for each $I_{t-1}^x \in \mathcal{I}_{t-1}^x(\Omega)$ and add these samples to $\tilde{\Omega}_t$. Larger values of $m$ will result in a larger set $\tilde{\Omega}_t$ from which we can form better value function estimates, but increase the time required to perform the backward pass. In this paper, we use $m=1$ and take one sample per information state at each time step in the backward pass. Note that using this method results in a constant sampling distribution for each iteration $k$ given by $\mathbb{Q}^k_t(w)=\frac{1}{|\mathcal{I}_{t-1}^x(\Omega)|}\sum\limits_{I_{t-1}^x \in \mathcal{I}_{t-1}^x(\Omega)} \mathbb{P}(w|I_{t-1}^x)$ for all $w \in \Omega_t$ and no updating step is necessary. We refer to this sampling method as \textit{standard sampling}.

An advantage to the standard sampling method is that it is very simple to implement. However, we encounter a problem in that high-risk, low-probability regions of $\Omega_t$ at each time step may not be sampled before the algorithm converges numerically. If the time horizon $T$ is quite large compared to the number of iterations necessary to reach convergence (as is the case in our application), then this is likely true for some values of $t$. If this occurs, the resulting policy may perform well for most of the sample paths, but when high-risk, low-probability outcomes do occur, the policy is likely to suffer. In the energy storage application, this would correspond to a situation where the policy does not plan enough storage to account for the possibility that wind will produce significantly less power than forecasted for an extended period of time, and, if this does occur, the system will experience a shortage.

This weakness of standard sampling motivates the need for \textit{risk-directed importance sampling}, a method that aims to sample these high-risk, low-probability regions of the outcome space with greater frequency by learning appropriate $\mathbb{Q}^k_t$'s. In our risk-directed importance sampling scheme, we seek to learn sampling distributions $\mathbb{Q}^k_t(w|I_{t-1}^x)$ for each $I_{t-1}^x \in \mathcal{I}_{t-1}^x(\Omega)$ that are ideally equal to $V_t^{x,*}(R_{t-1}^{x,k},w)\mathbb{P}(w|I_{t-1}^x)$ for $w \in \Omega_t$. As with standard sampling, to ensure $\mathbb{P}(w|I_{t-1}^x) >0$ for at least one element $w \in \tilde{\Omega}_t$, we sample $m$ elements (where we choose $m=1$ in this paper) from $\mathbb{Q}^k_t(w|I_{t-1}^x)$ for each $I_{t-1}^x \in \mathcal{I}_{t-1}^x(\Omega)$, where $\mathbb{P}(w|I_{t-1}^x) >0$ implies $\mathbb{Q}^k_t(w|I_{t-1}^x)>0$, and add these to $\tilde{\Omega}_t$. This implies that $\mathbb{Q}^k_t(w)=\frac{1}{|\mathcal{I}_{t-1}^x(\Omega)|}\sum\limits_{I_{t-1}^x \in \mathcal{I}_{t-1}^x(\Omega)} \mathbb{Q}^k_t(w|I_{t-1}^x)$ for all $w \in \Omega_t$. This technique therefore assigns a greater sampling probability to outcomes with large $V_t^{x,*}(R_{t-1}^{x,k},w)$ compared to the standard sampling technique. As risk is incorporated in the utility function by heavily penalizing cumulative shortages via the shortage state, this sampling method will thus visit risky regions of the outcome space more frequently.

Unfortunately, producing sampling distributions such that $\mathbb{Q}^k_t(w|I_{t-1}^x)= V_t^{x,*}(R_{t-1}^{x,k},w)\mathbb{P}(w|I_{t-1}^x)$ for $w \in \Omega_t$ is not a simple matter as: 1) $V_t^*(R_{t-1}^{x,k},w)$ is not known and is only approximated by $\ubar{V}^{x,k}_t(R_{t-1}^{x,k},w)$, 2) these approximations are both increasing and improving as the algorithm iteration count increases, and 3) we are unlikely to visit the same $R_{t-1}^{x,k}$ more than once in a high dimensional space. Since the value $\ubar{V}^{x,k}_t(R_{t-1}^{x,k},w)$ is dependent on $R_{t-1}^{x,k}$, we must account for this when learning sampling distributions. We therefore need an adaptive learning algorithm that overcomes these difficulties and continually updates our sampling distributions given the observations $\ubar{V}_t^{x,k}(R_{t-1}^{x,k},w)$ over the iterations of the algorithm. Our solution is based on the risk-directed sampling algorithm presented in \cite{jiang2017risk}.

We first address the effect of the resource state on the value $\ubar{V}^{x,k}_t(R_{t-1}^{x,k},w)$. Clearly this is dependent on $R_{t-1}^{B,x,k}$ as the cost of operating the system depends on how much energy is present in each storage device. Specifically, if the energy in the storage device is low, we are more likely to experience a shortage if a low wind power outcome $w \in \Omega_t$ occurs. Thus, we would like our sampling distributions to be dependent on $R_{t-1}^{B,x,k}$ as well, but encounter difficulties in that the resource state space is high-dimensional. To circumvent this issue, we consider an aggregated resource state based on the $\ell_1$-norm of the vector by partitioning $\mathbb{R}$ between $b_{t-1}^{0}=\min\limits_{R_{t-1}^{B,x}} ||R_{t-1}^{B,x}||_1$ and $b_{t-1}^{R-1}=\max\limits_{R_{t-1}^{B,x}} ||R_{t-1}^{B,x}||_1$ into $R$ bins at the division points $b_{t-1}^{0}<b_{t-1}^1<...<b_{t-1}^{R-1}$. The aggregate resource state $||R_{t-1}^{B,x}||_1$ belongs to bin $r$ if $b_{t-1}^{r} \leq ||R_{t-1}^{B,x,k}||_1< b_{t-1}^{r+1}$.

We then look to form the sampling distributions $\mathbb{Q}^k_t(w|I_{t-1}^x)$ using a weighted combination of basis distributions $\left\lbrace \phi_{trij} \right\rbrace_{j=0}^{J-1}$ for each time $t$, aggregated resource state ($||R_{t-1}^{B,x,k}||_1$) $r$, and information state $i \in \mathcal{I}^x_{t-1}(\Omega)$ as follows, 
\begin{align}
\mathbb{Q}^k_t(w|I_{t-1}^x=i)=\frac{1}{C^k_{tri}}\sum\limits_{j=0}^{J-1} \theta^k_{trij}(w) \phi_{trij}(w),
\end{align}where\begin{align}
C^k_{tri}=\sum\limits_{j=0}^{J-1} \sum\limits_{w \in \Omega_t} \theta^k_{trij}(w) \phi_{trij}(w).
\end{align}
We require that $\sum\limits_{w \in \Omega_t} \phi_{trij}(w)=1$, $\phi_{trij}(w)>0$ if and only if $\mathbb{P}(w|I_{t-1}^x=i)>0$, and $\phi_{tri0}=\mathbb{P}(w|I_{t-1}^x=i)$ for all $t$ and $i \in \mathcal{I}^x_{t-1}(\Omega)$. The vector of weights at iteration $k$, $\theta^k_{tri}$, are ideally chosen using the projection operator $\Pi_{\theta}[|V_t^{x,*}(R_{t-1}^{x,k},w)|\mathbb{P}(w|I_{t-1}^x=i)]$, where, letting $W^u_t$ be a uniformly distributed random variable over the outcome space $\Omega_t$,
\begin{align}
\Pi_{\theta}F=\argmin\limits_{\theta \geq 0}\mathbb{E}\left[\left[\theta^T \phi(W^u_t)-F(W^u_t)\right]^2\right].
\end{align}
If the weights are chosen in this manner, $(\theta^k_{tri})^T \phi(w)$ achieves the best fit (in terms of minimizing the expected mean squared error) to $|V_t^{x,*}(R_{t-1}^{x,k},w)|\mathbb{P}(w|I_{t-1}^x)$ for $w \in \Omega_t$ given the set of basis distributions.

However, we do not know $V_t^{x,*}(R_{t-1}^{x,k},w)$ and can only observe the lower bound approximation $\ubar{V}_t^{x,k}(R_{t-1}^{x,k},w)$ at each iteration $k$. Note that $\bar{V}_t^{x,k}(R_{t-1}^x,I_t^x=i)$ (and thus $\ubar{V}^{x,k}_t(R_{t-1}^{x,k},w)$) may increase greatly as $k$ increases, especially in the first few iterations. Therefore, it is entirely possible that we have $w_m \in \Omega_t$ sampled at iteration $k'$ and $w_n \in \Omega_t$ sampled at iteration $k''>k'$, such that for the same point $R_{t-1}^{x,k'}=R_{t-1}^{x,k''}$, $V_t^{x,*}(R_{t-1}^{x,k'},w_m)>V_t^{x,*}(R_{t-1}^{x,k''},w_n)$, but $\ubar{V}^{x,k'}_t(R_{t-1}^{x,k'},w_m)<\ubar{V}^{x,k''}_t(R_{t-1}^{x,k''},w_n)$. Since our goal is to sample outcomes $w \in \Omega_t$ that result in higher values of $V_t^{x,*}(R_{t-1}^{x,k},w)$ with more frequency than their nominal probability, using these estimates would produce undesirable sampling distributions for subsequent iterations. For this reason, we define, for $w \in \Omega_t$,
\begin{align}
v_t^{x,*}(R_{t-1}^x,w)=V_t^{x,*}(R_{t-1}^x,w)-\sum\limits_{i \in \mathcal{I}_t^x(\Omega)} V_t^{x,*}(R_{t-1}^x,I_t^x=i)\mathbb{P}(I_t^x=i|w),
\end{align} and an observation of $v_t^{x,*}(R_{t-1}^x,w)$ at iteration $k$,
\begin{align}
\hat{v}_t^{x,k}(R_{t-1}^x,w)= \max \left\lbrace 0,\ubar{V}_t^{x,k} (R_{t-1}^x,w)-\sum\limits_{i \in \mathcal{I}_t^x(\Omega)} \bar{V}_t^{x,k}(R_{t-1}^x,I_t^x=i)\mathbb{P}(I_t^x=i|w) \right\rbrace,
\end{align} which, for a realization of $w \in \Omega_t$, estimates the one-step cost plus the \textit{change} in the expected downstream value as a result of altering the resource state. Note that errors in the current VFA may cause the expression to dip below zero in certain cases, thus we enforce that $\hat{v}_t^{x,k}(R_{t-1}^x,w)$ be nonnegative. This variable better isolates the impact of the outcome $w \in \Omega_t$ on the system with less dependence on the iteration count or the accuracy of the VFAs. We thus aim to fit sampling distributions to $v_t^{x,*}(R_{t-1}^{x,k},w)\mathbb{P}(w|I_{t-1}^x=i)$ and choose weights using the projection $\Pi_{\theta}[|v_t^{x,*}(R_{t-1}^{x,k},w)|\mathbb{P}(w|I_{t-1}^x=i)]$.

Though improved, $\hat{v}_t^{x,k}(R_{t-1}^x,w)$ is still dependent on $k$, and we are only able to obtain observations for relatively few outcomes $w \in \Omega_t$ at each iteration $k$ (sampling too many outcomes such that $\tilde{\Omega}_t$ is close to $\Omega_t$ in size would defeat the purpose of sampling). However, if the weights $\theta^k_{tri}$ are refit using only the observations in the current iteration, this may produce vastly different sampling distributions from iteration to iteration. Thus, to update the weights based on new observations at iteration $k$ for $i \in \mathcal{I}^x_{t-1}(\Omega)$, we utilize the weights from the previous iteration and the batch recursive updating formula $\theta^k_{tri}=\max\left\lbrace \bf{0},\theta^{k-1}_{tri}-\gamma_{tr}^k (\Phi^{k}_{tri})^T\varepsilon_{tri}^k\right\rbrace$, where, if $L$ samples were taken in total at time $t$ during iteration $k$, \begin{align*}
\Phi^{k}_{tri}=\begin{bmatrix}
\phi_0(w_0) & \phi_1(w_0) & \ldots & \phi_{J-1}(w_0) \\
\phi_0(w_1) & \phi_1(w_1) & \ldots & \phi_{J-1}(w_1) \\
\vdots& \vdots&\ddots&\vdots\\
\phi_0(w_{L-1}) & \phi_1(w_{L-1}) & \ldots & \phi_{J-1}(w_{L-1})
\end{bmatrix},
\end{align*} and the $\ell$-th element of the length $L$ column vector $\varepsilon_{tri}^k$ is given by \begin{align*}
\varepsilon_{tri}^k(\ell)=\frac{(\theta^k_{tri})^T\phi_{tri}(w_{\ell})-|\hat{v}_t^{x,k}(R_{t-1}^{x,k},w_{\ell})|\mathbb{P}(w_{\ell}|i)}{\mathbb{Q}^{k-1}_t(w_{\ell})}.
\end{align*} In the special case that $\theta^k_{tri}=\bf{0}$, all basis distributions are weighted equally.

This mimics the form of the updating formula for recursive ordinary least squares regression, but allows us to control the stepsize rule. We choose a stepsize of $\gamma_{tr}^k=\frac{a}{a+n^k}$, where $n^k$ is the count of how many times the aggregate resource state $r$ has been visited at time $t$, and $a$ is a constant that is tuned to work well for the problem (note that more sophisticated stepsize rules can be used instead). This is utilized instead of the common $\gamma_{tr}^k=\frac{1}{n^k}$ stepsize as it places more weight on recent observations of $v_t^{x,k}(R_{t-1}^x,w)$.

\begin{figure}
\centering
\includegraphics[width=\columnwidth, height= 4.5 in]{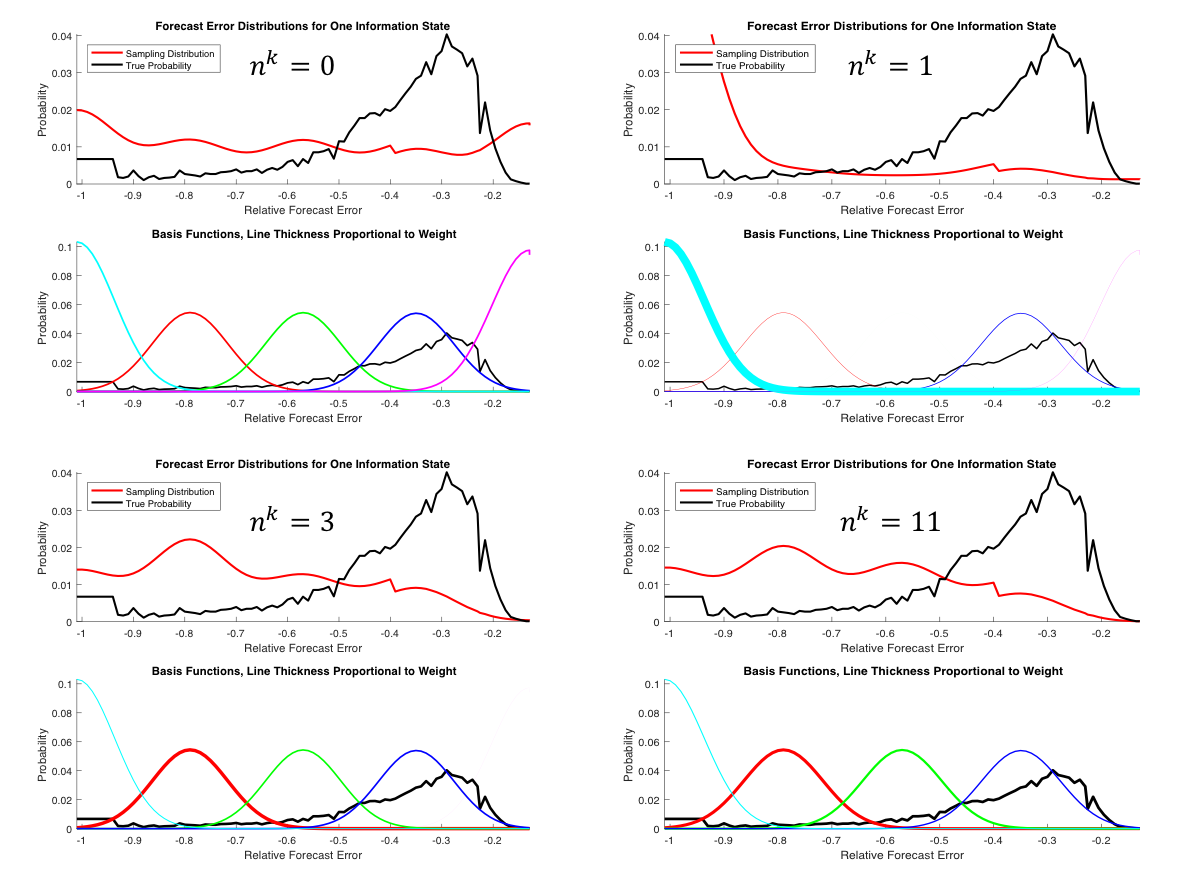}
\caption{We show the evolution of a sampling distribution for a single time $t$, aggregated resource state $r$, and information state $i$. $n^k$ represents the number of times this state has been visited by the algorithm. The top graph in each quadrant shows the true distribution $\mathbb{P}(w|I_{t-1}^x=i)$ versus the conditional sampling distribution $\mathbb{Q}_t^k(w|I_{t-1}^x=i)$ for $w \in \Omega_t$ where each outcome is a wind forecast error (as a fraction of the maximum output of the wind source). The bottom graph in each quadrant shows the basis distributions, where line weight is proportional to the weight placed on each distribution. Note that as the iteration count increases, we learn to sample the larger magnitude negative errors with more frequency from as these are more likely to produce shortages.}
\label{IS_Sampling_Dists}
\end{figure}

In this paper, for each time $t$, aggregated resource state $r$, and information state $i$, the basis distributions $\left\lbrace \phi_{trij} \right\rbrace_{j=0}^{J-1}$ are chosen to be a set of normal and half-normal distributions with different means and variances that span the support of $\mathbb{P}(W_t|I_{t-1}^x=i)$ along with $\mathbb{P}(W_t|I_{t-1}^x=i)$ itself. This choice of basis distributions generalizes quite well to fit a wide variety of arbitrary shapes for $|v_t^{x,*}(R_{t-1}^{x,k},w)|\mathbb{P}(w|I_{t-1}^x)$ for $w \in \Omega_t$ when properly weighted. Figure \ref{IS_Sampling_Dists} shows an example of the evolution of a sample distribution over iterations of the algorithm for a single time $t$, aggregated resource state $r$, and information state $i$. The top graph in each quadrant shows the true distribution $\mathbb{P}(w|I_{t-1}^x=i)$ versus the conditional sampling distribution $\mathbb{Q}_t^k(w|I_{t-1}^x=i)$ for $w \in \Omega_t$. The bottom graph in each quadrant shows the set of basis distributions (an intentionally small set for illustration purposes) where line thickness is proportional to the weight placed on each distribution. Note that as the algorithm progresses, we learn to sample with more frequency from the larger magnitude negative forecast errors which produce higher values $\hat{v}_t^{x,k}(R_{t-1}^{x,k},w)$ as they are more likely to lead to shortages.

\section{Numerical Results}
\label{Numerical Results}
In this section we present results across a variety of test cases that support the following claims:
\begin{henumerate}

\item Power grids with high penetrations of renewables can reduce shortages by intelligently operating a system of distributed battery storage devices to smooth out the variability present in renewables.

\item Shortages can be further reduced if renewables are modeled with a crossing state model which captures the distribution of times for which the stochastic process is above or below its forecast. This is in comparison to a stochastic model which assumes intertemporal independence.

\item Since the SDDP algorithm must be modified to produce a set of Bender's cuts for each possible information state at each time step to accommodate the unique stochastic model, the algorithm is now too computationally inefficient to be used in practice. A simple sampling scheme can be employed in the backward pass to speed up the convergence of the SDDP algorithm, but doing so results in diminished solution quality. In particular, it results in more shortages than when the full outcome space is considered at each time step. Using an importance sampling method, we achieve accelerated convergence without sacrificing solution quality or robustness.

\end{henumerate}

\begin{figure}
\centering
\includegraphics[width=\columnwidth, height=3.25 in]{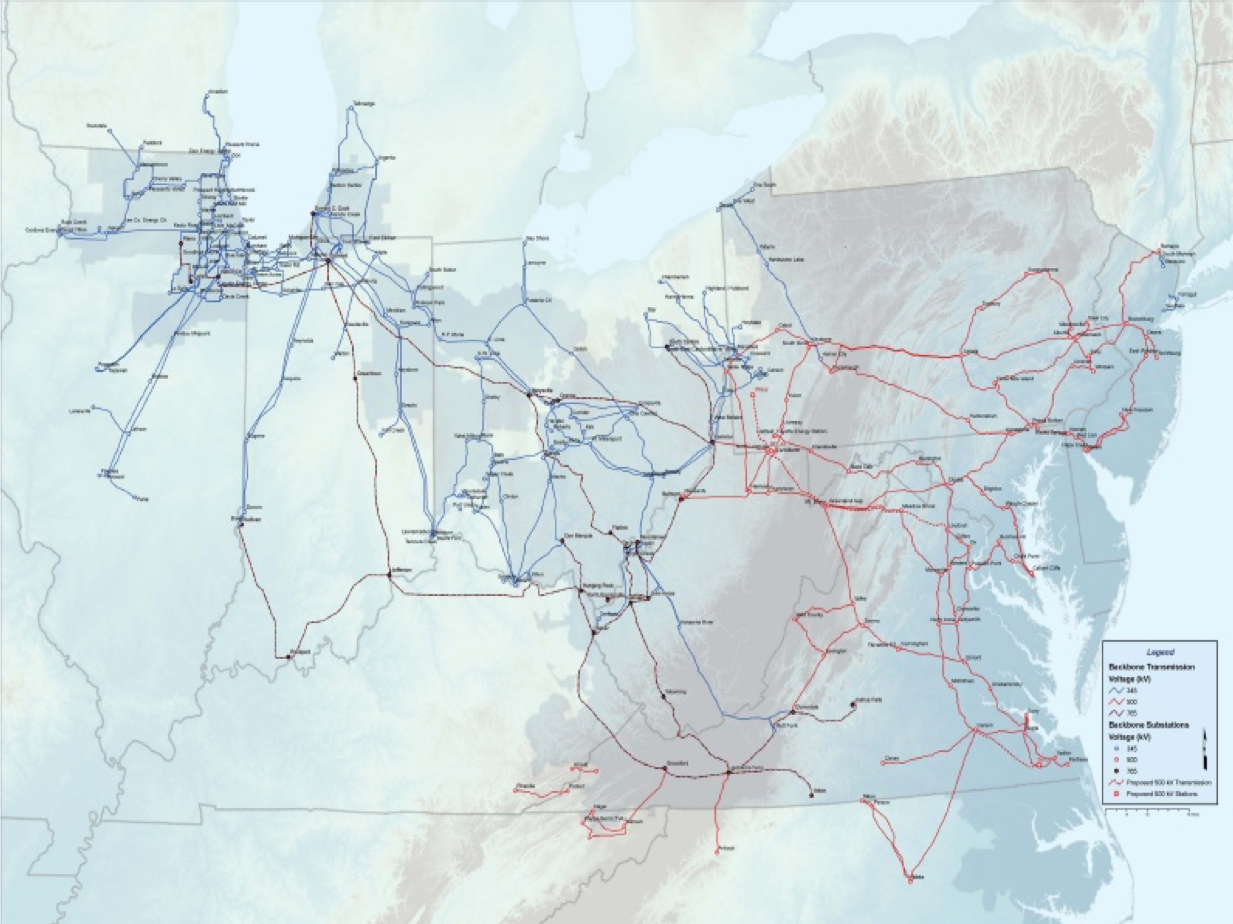}
\caption{The geographical area covered by PJM. Also shown are the largest capacity transmission lines in PJM.}
\label{PJM}
\end{figure}

The system used to test these claims is a carefully calibrated network model of the PJM grid called SMART-Storage that simulates grid operation while optimizing distributed storage. For simplicity, we restrict our study to the highest voltage lines, which forms an aggregated grid model with 1,360 buses and 1,715 transmission lines (instead of the full 9,000 buses and 14,000 lines). A map of the geographical area covered by PJM, along with the very largest capacity PJM transmission lines is shown in Figure \ref{PJM}. The generators that make up the set $\mathcal{G}$ consist of 396 gas turbines with a combined maximum power capacity of 23,309 MW, 50 combined cycle generators with a capacity of 21,248 MW, 264 steam generators with a capacity of 73,374 MW, 31 nuclear reactors with a capacity of 31,086 MW, and 84 conventional hydro power generators with a power capacity of 2,217 MW. Using a demand curve taken from a typical PJM day we experience an average demand of 91,054 MW, morning and afternoon peaks with maximums of 101,743 MW and 98,695 MW respectively, and nighttime and midday troughs with minimums of 80,902 MW and 86,454 MW respectively.

SMART-Storage requires that the discrete on-off unit commitment decisions for generators over the optimization horizon are input to the system. These are the $Z_{t,g}^G$ variables defined in Section \ref{StateVar}. For this, SMART-Storage leverages a previous model developed by \cite{simao2017challenge} called SMART-ISO that can simulate day-ahead and real-time planning processes at PJM, mimicking actual operations. The model has been calibrated using historical performance at PJM and was shown to accurately replicate network behavior. SMART-ISO is not designed to include and optimize distributed storage on a large scale, but is useful for its high quality unit commitment decisions.

Prior to delving into the details of the studies, we note that the system has undergone extensive tuning of many underlying parameters to produce numerically useful results. These include inputs such as the relative primal feasibility tolerance and relative complementary gap tolerance of the convex optimization solver (IBM ILOG CPLEX 12.4) used to find solutions to each subproblem. It is important to note that inappropraite choices for these tolerances may produce stopping conditions in the solver which cannot be satisfied. A thorough discussion of how to choose these tolerances when using a regularized SDDP with Markov uncertainty algorithm for a problem with long time horizons is available in \cite{asamov2015regularizedARXIV}. The regularization sequence $\left\lbrace \rho^k \right\rbrace$ (see Section \ref{StagewiseDependentSDDP}) also has an effect on the solution and is chosen to be the geometric sequence $\rho^k=\rho^0 r^k$ where $r=0.95$ and $\rho^0=1$ as suggested by \cite{asamov2015regularized}. Sensitivity analysis on these parameters is not within the scope of this paper, however \cite{asamov2015regularized} presents empirical results showing that several different pairs of $(r,\rho^0)$ with $r$ ranging from $.9$ to $.99$ and $\rho^0$ from $1$ to $100$ produce similar convergence behavior. Once set, these underlying parameters are held constant in our studies, regardless of the model or sampling method chosen.

\subsection{The Value of Battery Storage and the Crossing State Stochastic Model for Renewables}
\label{Value of HSMMs}

To simulate and observe the effects of increased penetrations of renewables, we introduce offshore wind into the system with 16,189 MW of power capacity, which represents a significant portion of the instantaneous demand at any point during the day. Daily day-ahead offshore wind forecasts are generated using a meteorological model known as the Weather, Research, and Forecasting (WRF) model. Wind forecast error models are trained using sets of forecasted and corresponding actual time series of wind power collected from a set of PJM wind farms over different months of the 2013. The total wind capacity is then adjusted to the desired the maximum power output by scaling both the errors and forecasts proportionally. For a more detailed description and validation of this procedure see \cite{archer2017challenge}, as we utilize the same method but assume different stochastic wind models. We discretize the wind forecast error outcome space $\Omega_t$ into roughly $200$ outcomes per time step depending on the forecast level $f_t^E$. Periods for which the forecast is near the minimum/maximum power level will have fewer negative/positive possible forecast error outcomes.

We also place 20 large storage devices across the grid at points of interest, such as at points where grid congestion is often experienced without storage and near points of interconnection for offshore wind. These have varying capacities, efficiencies, and charge/discharge rates but are all large enough to have a significant impact on grid operations. In each test case, we develop the policy to control these batteries by training VFAs using the appropriate SDDP algorithm until the empirical upper bound and lower bound have converged to within two percent of each other. Then, 50 testing iterations are carried out in which we fix the value functions and observe the performance of the policy over 50 different wind forecast error sample paths.

\begin{figure}
\centering
\includegraphics[width=\columnwidth, height= 2.4 in]{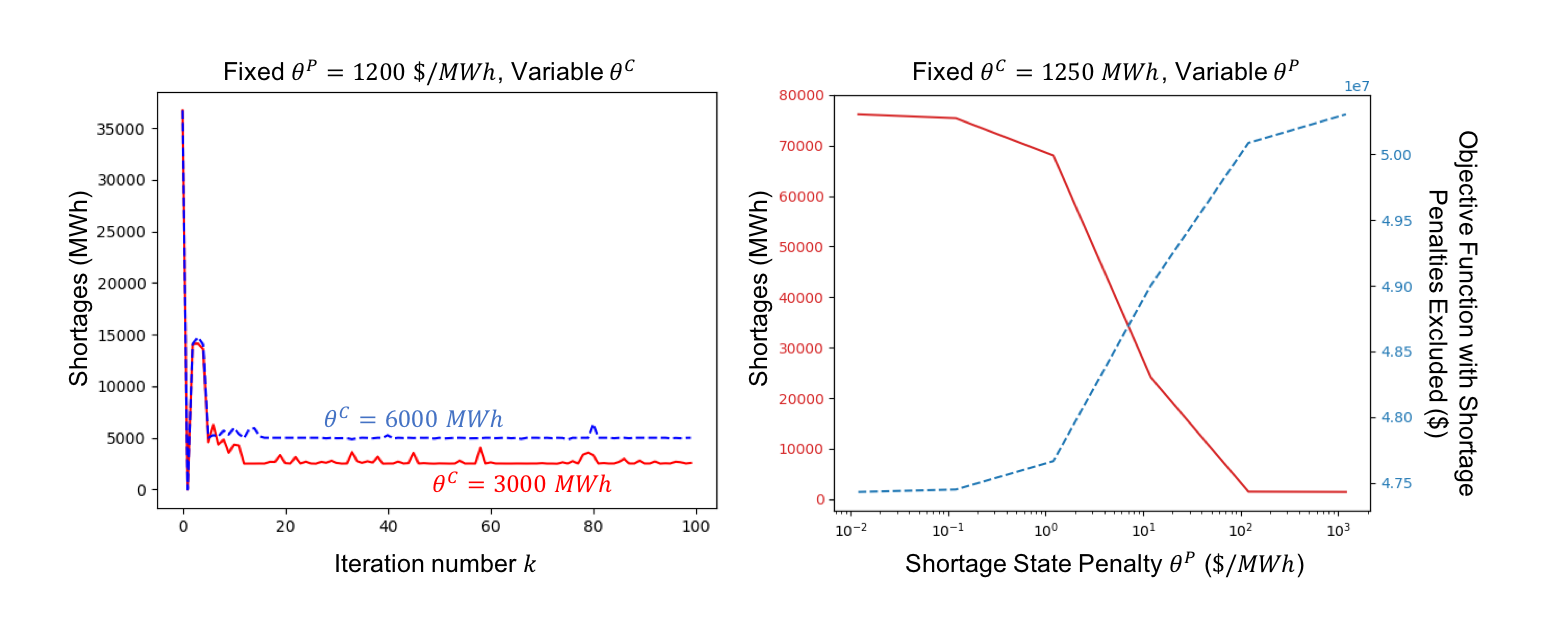}
\caption{Left: the number of shortages over 100 iterations of the algorithm. Observe that varying the shortage state threshold will directly affect the number of shortages observed when simulating the policy. Provided that $\theta^P$ is large enough to discourage shortages in excess of the $\theta^C$, we see convergence to a solution in which cumulative shortages are below $\theta^C$ in most cases. Right: Shortage state penalty $\theta^P$ versus average shortages and objective function cost (minus shortage penalties) over 50 testing iterations after convergence. Here we see that finding an appropriate $\theta^P$ for a given $\theta^C$ is crucial as well. If $\theta^P$ is too small, the threshold is essentially ignored, while a large $\theta^P$ may drive up generation costs without producing a solution that is noticeably more shortage-averse.}
\label{Cost_v_Shortage_Curve}
\end{figure}

While most of the remaining parameters, such as the grid structure, power capacities of generators, and costs for each form of generation are fixed and come from the SMART-Storage model of the PJM grid, a few tunable parameters can have a great effect on the behavior and performance of the system. Most notably for this study, the shortage state threshold $\theta^C$, and the associated penalty $\theta^P$ placed on shortages in excess of the threshold, will have large impacts on the risk-cost tradeoff in the system. In the left plot of Figure \ref{Cost_v_Shortage_Curve}, we see that the solution (computed by Algorithm \ref{alg:SDDP} with importance sampling) converges to a point such that $R_T^S<=\theta^C$ in most scenarios provided that $\theta^P$ is large enough. In the right plot of Figure \ref{Cost_v_Shortage_Curve}, we observe that finding the appropriate $\theta^P$ to accompany $\theta^C$ is important as well. If $\theta^P$ is too small, the threshold will be largely ignored as other components of the objective function will outweigh the shortage state penalty. Conversely, a much too large $\theta^P$ will produce a policy that may be too conservative in some cases. For example, in Figure \ref{Cost_v_Shortage_Curve}, we incur higher generation costs when increasing $\theta^P$ from $\$$120/MWh to $\$$1200/MWh without producing a solution that results in noticeably fewer shortages.

Though accounting for shortage risk in this manner (i.e. by augmenting the state variable to track shortages and applying a threshold penalty at time $T$) does produce the need to tune these parameters, it is actually quite practical from a system operator's point of view as the system's aversion to risk can be directly and intuitively controlled. Let us assume the role of a grid operator and decide that we want to operate with a low tolerance for shortages. Thus, we set a high penalty $\theta^P=960$ $\$$/MWh and test the system under two relatively small thresholds: $\theta^C=10$ MWh and $\theta^C=500$ MWh. While further sensitivity analysis could be performed on these (and other) system parameters, this is not a focus of the paper. We instead now show that solution quality and robustness can be improved by choice of stochastic model and SDDP algorithm.

The upper half of Table \ref{table:model and sampling comparison} (separated by the double horizontal line) displays statistics for the number of shortages and the objective function over testing iterations under three different modeling assumptions:
\begin{henumerate}
\item The power system run without any storage devices. In this case no SDDP algorithm is necessary as there is no controllable resource state.
\item With distributed storage devices and the exogenous wind process modeled with an IID model. In this case standard SDDP with quadratic regularization is utilized to fit value functions.
\item With distributed storage devices and the exogenous wind process modeled with the univariate crossing state model from \cite{durante2017b}. In this case the version of SDDP from Algorithm \ref{alg:SDDP} \textit{without} sampling is utilized (i.e. $\tilde{\Omega}_t=\Omega_t$ and $\mathcal{L}(w|I_{t-1}^x)=1$ for $w \in \Omega_t$). We aim to compare stochastic models here, and will examine the impacts of sampling and choice of sampling method in Subsection \ref{IS Numerical}.
\end{henumerate} To account for the fact that performance is dependent on the scenarios used for both training and testing, for each $\theta^C$ we repeat the following training and testing procedure over four different random seeds: 1) 50 wind scenarios are generated and stored for testing. 2) For each model, the policy is trained until convergence, at which point the VFAs are fixed. The time required to reach numerical convergence is recorded in seconds. 3) Using the set of test scenarios (which are consistent across the models), 50 test iterations are carried out for each model. We gather the mean, standard deviation (SD), and worst case scenario (WC) of both the objective function and cumulative shortages experienced by the system per iteration. Table \ref{table:model and sampling comparison} reports the \textit{average} value of each statistic over all the repetitions of the training/testing procedure for each model and $\theta^C$ combination.

To properly compare the performance of the models, one should use sets of results with the same $\theta^C$. To help the reader differentiate, the results for $\theta^C=500$ MWh are italicized in Table \ref{table:model and sampling comparison}, while those for $\theta^C=10$ MWh are not. When results are compared in this manner, Table \ref{table:model and sampling comparison} clearly shows that intelligently controlled distributed storage devices can drastically reduce shortages in a power system with high penetrations of renewables as both models that utilize storage vastly outperform the model without storage devices. We also see that, while the IID modeling assumption may perform similarly on average in terms of minimizing the objective function, it is less robust to scenarios when wind underperforms its forecast and suffers more shortages than the crossing state wind model. Observe that the worst case scenario for both the objective function and shortages are also less extreme for the crossing state model, as the policy better prepares for the possibility of poor wind scenarios.

\begin{table}
\TABLE
{Comparison of result statistics over several test cases for the various models and algorithms discussed in this paper. The testing and training procedure described in Section \ref{Value of HSMMs} is carried out for different values of $\theta^C$ and the \textit{average} value of each statistic (where SD = standard deviation and WC = worst case) over all the training/testing cases is reported. One should compare the performance of the models by comparing sets of results with the same $\theta^C$. To help differentiate, the results for $\theta^C=500$ MWh are italicized, while those for $\theta^C=10$ MWh are not. \label{table:model and sampling comparison}}
{\begin{tabular}{c || c | c | c c c | c c c } 
\hline Model/ & \multirow{2}{*}{\parbox{.8cm}{$\theta^C$ in MWh}} & \multirow{2}{*}{\parbox{.8cm}{Covg. Time}} & \multicolumn{3}{c}{Objective Function Stats ($\$$)} \vline & \multicolumn{3}{c}{Shortage Stats (MWh)} \\
\cline{4-9}
Method & & & Mean ($\times 10^7$) & SD ($\times 10^5$) & WC ($\times 10^7$) & Mean & SD & WC\\
\hline
No Storage/ & 10 & N/A & 6.7372 & 5.8082 & 6.9713 & 41459.17 & 960.80 & 45257.12\\
\cline{2-9} N/A & \textit{500} & \textit{N/A} & \textit{6.7395} & \textit{5.6985} & \textit{6.9765} & \textit{41471.61} & \textit{935.53} & \textit{45290.89}\\
\hline
IID Model/ & 10 &  4491 & 6.0026 & 2.7639 & 6.1380 & 137.69 & 202.16 & 1240.25\\
\cline{2-9} No Sampling & \textit{500} & \textit{ 5216} & \textit{5.9825} & \textit{1.5669} & \textit{6.0614} & \textit{536.34} & \textit{86.63} & \textit{1027.50}\\
\hline
Crossing State/ & 10 &  25708 & 5.9896 & 1.9597 & 6.0842 & 43.18 & 109.39 & 628.83\\
\cline{2-9} No Sampling & \textit{500} & \textit{ 30851} & \textit{5.9750} & \textit{0.9854} & \textit{6.0253} & \textit{508.31} & \textit{31.03} & \textit{662.67}\\
\hline
\hline
Crossing State/ & 10 &  8385 & 5.9968 & 2.1688 & 6.0939 & 43.01 & 118.98 & 647.23\\
\cline{2-9} Standard Sampling & \textit{500} & \textit{ 9602} & \textit{5.9796} & \textit{1.1000} & \textit{6.0257} & \textit{534.58} & \textit{47.83} & \textit{726.93}\\
\hline
Crossing State/ & 10 &  5938 & 6.0007 & 1.6824 & 6.0669 & 26.21 & 54.06 & 322.97\\
\cline{2-9} Importance Sampling & \textit{500} & \textit{ 7246} & \textit{5.9893} & \textit{1.0393} & \textit{6.0219} & \textit{525.40} & \textit{17.49} & \textit{553.61}\\
\hline

\end{tabular}}
{}
\end{table}

\subsection{Accelerating the Convergence of SDDP with Markov Uncertainty through Importance Sampling}
\label{IS Numerical}

One downside of using the crossing state stochastic model and associated SDDP algorithm is the increased computation time required to reach numerical convergence. Two main factors account for this. In comparison to the IID model: 1) we must now fit value functions $\bar{V}_t^{x,k}(R_t^x,I_t^x)$ for $I_{t}^x \in \mathcal{I}_{t}^x(\Omega)$ instead of one VFA for the resource state and 2) in order to realize the benefits of a more sophisticated model, the outcome space $\Omega_t$ must be discretized more finely than when using the IID model. To understand why the latter point is true, consider the difference between discretizing $|\Omega_t|$ into 15 or 200 outcomes. As long as the points are properly chosen (say evenly on both sides of the wind forecast), the accuracy of the observation of the value function of the resource state ($\ubar{V}_{t+1}^{x,k}(R_t^x)$ for the IID model) will only marginally improve with larger $|\Omega_t|$. However, with small $|\Omega_t|$, we lose key information about the true dynamics of the stochastic process (stored in $I_t^x$), and thus lose the value of using a more sophisticated stochastic model (we note that the model comparison in Section \ref{Value of HSMMs} was, however, performed on the same discretization of $\Omega_t$, this explanation is solely to point out that the IID model performance does not benefit as much from larger $|\Omega_t|$). A larger $|\Omega_t|$ requires solving more linear programs per time step in the backward pass, one for each $w \in \Omega_t$, which is the computational bottleneck in our algorithm. Computation time constraints for the system may restrict the use of the more sophisticated model. We thus present results on how intelligently sampling $\Omega_t$ can decrease computation time without sacrificing solution quality.

\begin{figure}
\centering
\includegraphics[width=\columnwidth, height= 2.3 in]{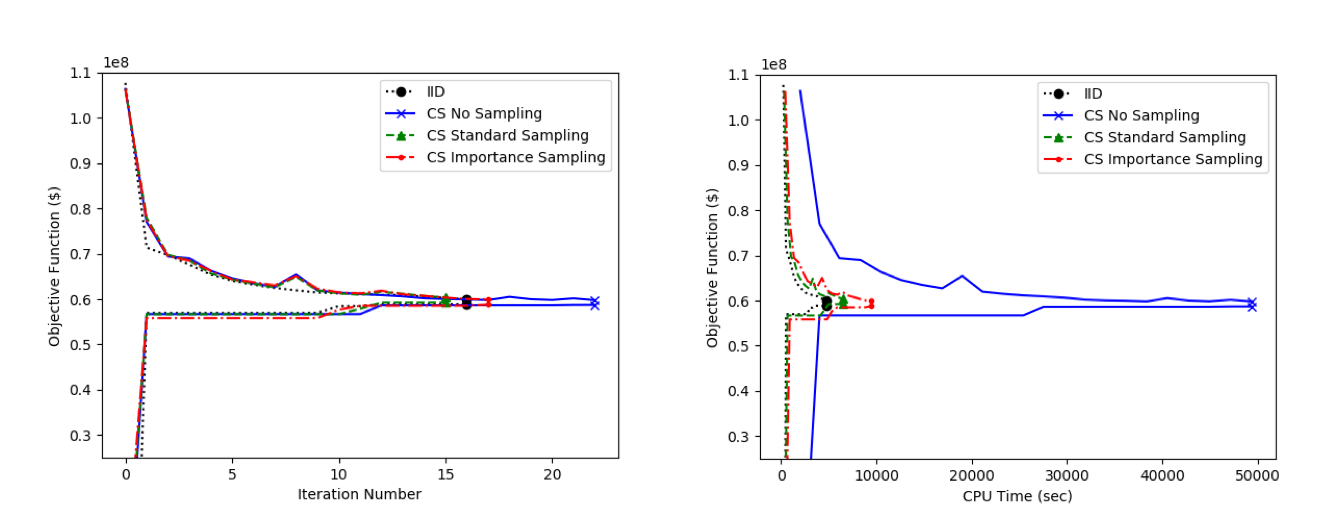}
\caption{Left: For the selected models and sampling types, the lower bound cost and the estimate of the upper bound cost from the forward pass are plotted versus the iteration number in one test case. Convergence is reached when the bounds are within $\epsilon=2$ percent of each other. Right: The same set of training iterations are shown, except the evolution of the bounds are plotted with respect to CPU time here. Note that while the number of iterations to converge may be roughly the same, the total CPU time required to use the crossing state model is drastically reduced via sampling. The sampling algorithms are nearly on par with the IID model in terms of computation time required for numerical convergence.}
\label{convergencePlot}
\end{figure}

We extend the studies performed in Section \ref{Value of HSMMs}, by adding results obtained by sampling the outcome space in the backward pass. Two sets of additional results are shown in the lower half of Table \ref{table:model and sampling comparison} (below the double horizontal line). The crossing state wind model and the version of SDDP from Algorithm \ref{alg:SDDP} are used in each case, but the sampling method is varied. The ``Standard Sampling" and ``Importance Sampling" methods correspond to the methods described in Section \ref{IS Section}. These both sample approximately $15$ percent of $\Omega_t$ per time step. As a result, using the numbers from Table \ref{table:model and sampling comparison}, we see that numerical convergence is accelerated by factors of $3.07$ and $3.21$ with standard sampling and $4.33$ and $4.26$ with importance sampling. Figure \ref{convergencePlot} displays the convergence of the solution across different models for one test case graphically. Observe that the number of iterations until convergence (left plot) is not altered much by sampling, but the total CPU time (right plot) is reduced as each iteration now takes less time to complete. Additionally, the sampling algorithms are nearly as efficient as the IID model in terms of computation time required for numerical convergence.

Though both sampling methods achieve the goal of expediting convergence, importance sampling produces higher quality solutions. This is supported by the results presented in Table \ref{table:model and sampling comparison} where, again, we compare results from each $\theta^C$ separately. Note that the overall performance of the importance sampling method is comparable to that of the unsampled crossing state model, but there are several key differences in the results. While the average objective function values are slightly higher, the standard deviations and worst case scenarios are lower. This indicates that the policy produced with importance sampling is more robust to poor wind power scenarios. This is supported by the fact that the importance sampling method results in fewer shortages, both on average and in the worst case. Furthermore, when these same comparisons are made to the IID model from Section \ref{Value of HSMMs}, the ability of the crossing state model with importance sampling to produce risk-averse policies is even more pronounced.

\begin{figure}
\centering
\includegraphics[width=\columnwidth, height= 5 in]{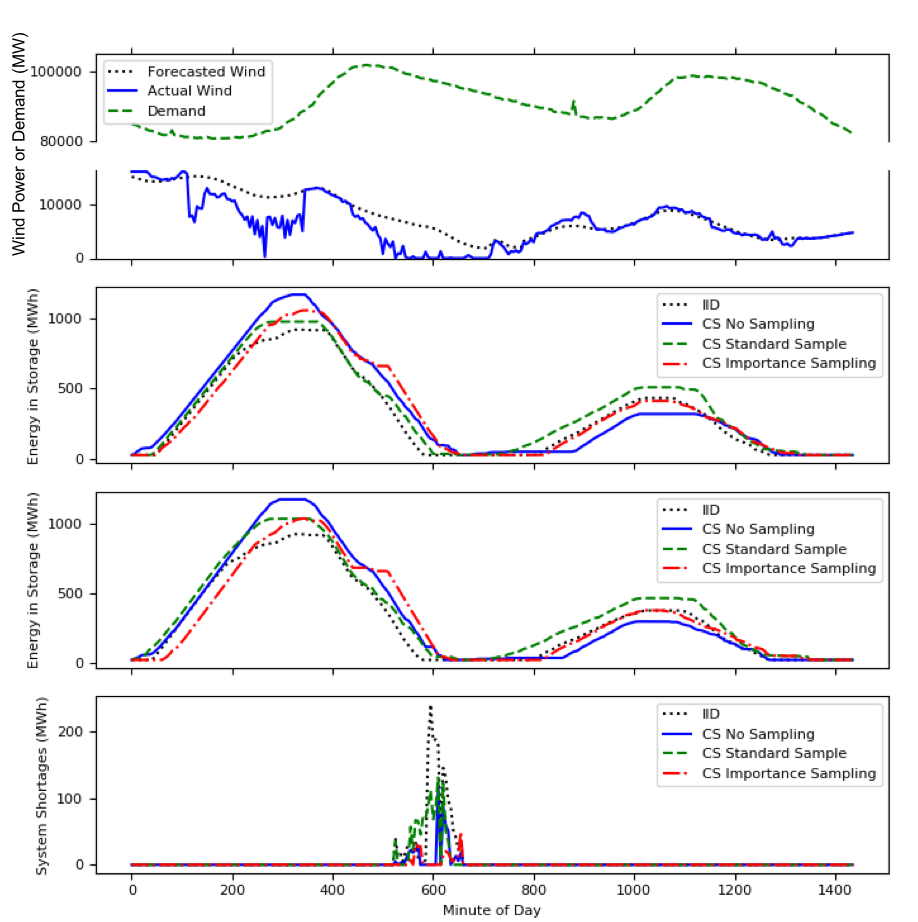}
\caption{Top: The demand, the forecasted wind power, and one sample path of actual wind power output over one day (note the break in the y-axis). Middle Plots: The energy in two representative batteries following the policy produced by each model and algorithm under the wind scenario shown in the top graph. Bottom: The total system shortages experienced when following the policy produced by each model and algorithm under the wind scenario shown in the top graph. Observe that the policy developed using the IID model has less power in storage during critical periods, as does the standard sampling version of the crossing state model (to a lesser extent).  Conversely, the unsampled crossing state model, and especially the importance sampled version, will plan extra storage to prepare for the possibility that wind may drops below its forecast for an extended period of time. Consequently, when such a scenario occurs, we see the policies developed with the IID model and crossing state model with standard sampling experiencing more shortages as they run out of backup storage more quickly.}
\label{DemandWindStorageShortage}
\end{figure}

Conversely, the quality of the policy produced by the standard sampling method is poor in comparison as it is both higher in mean objective value and less risk-averse than the unsampled version. The standard sampling method is heavily dependent on which outcomes $w \in \Omega_t$ are sampled in the training iterations. While the importance sampling method seeks to sample outcomes that produce high values of $\hat{v}_t^{x,k}(R_{t-1}^x,w)$ with higher probability, standard sampling does not sample intelligently. If high-risk outcomes are not sampled during training, but then appear in practice (testing), the policy may not have stored enough energy to prevent shortages. Figure \ref{DemandWindStorageShortage} shows, for one testing iteration, a plot of energy storage and shortages over the course of one day in five-minute time intervals for each sampling method highlighting such an occurrence (the IID model is also included and exhibits similar behavior to the standard sampling method). The result is a less risk-averse solution which, when compared to the unsampled and importance sampling methods, produces high worst case objective function values and cumulative shortages that often far exceed the threshold set by the system operator $\theta^C$ (see Table \ref{table:model and sampling comparison}).

\begin{figure}
\centering
\includegraphics[width=\columnwidth, height= 2.3 in]{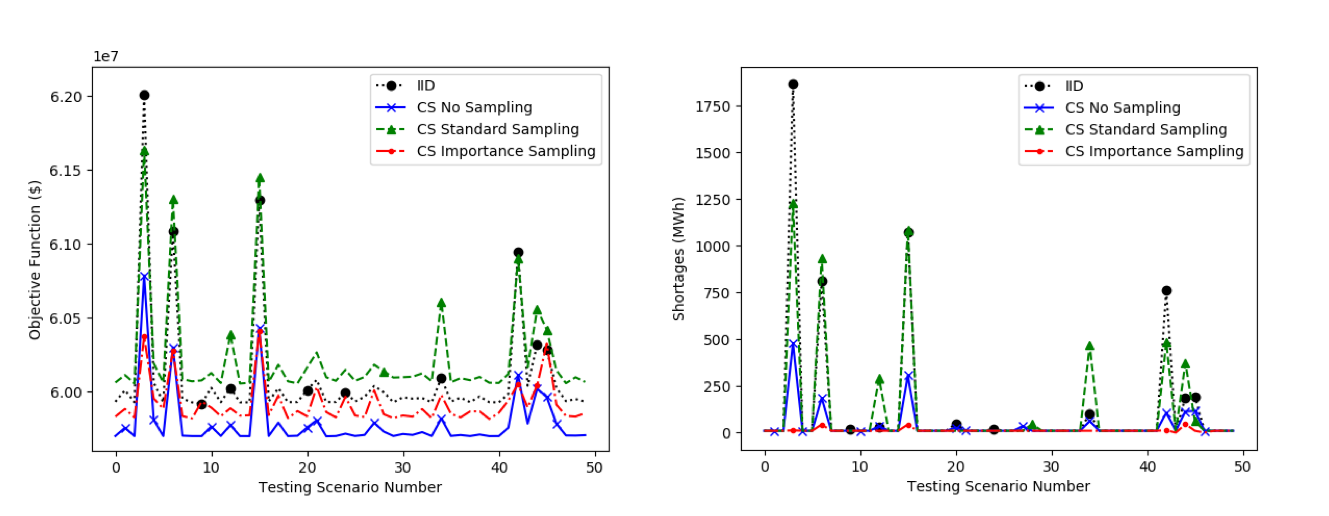}
\caption{Left: For the selected models and sampling types, the objective function is plotted over 50 testing scenarios in one test case. Right: For the same set of scenarios, the cumulative shortages are plotted in MWh. In both plots, a marker is placed when $R_T^S>\theta^C$ for a given scenario. We observe that the crossing state model performs better than the IID model in terms of objective function and aversion to shortages. Applying standard sampling may speed up convergence but reduces solution quality. In this case, the risk-averse importance sampling algorithm actually results in fewer shortages compared to the unsampled version, but does perform slightly worse in terms of overall cost-efficiency.}
\label{ObjAndShortVModel}
\end{figure}

The claims made about the quality of each solution method discussed in the previous two paragraphs are reinforced by the plots in Figure \ref{ObjAndShortVModel}, which displays both the objective function and cumulative shortages for one test case with $\theta^C=10$ MWh over 50 iterations. Here we can see graphically that the objective function performance is lowest on average for the unsampled crossing state model, but the importance sampling method results in the fewest shortages and lowest worst case objective function. Furthermore, both the IID model and standard sampling method both experience scenarios in which not enough storage was planned and cumulative shortages greatly exceed $\theta^C$. Additionally, they produce higher mean and worst case objective values.

\section{Conclusions}
\label{Conclusion}

Grid-level distributed storage can help mitigate risk in power systems with significant penetrations of renewable energy. To do so most effectively, the system control policy should accurately model and prepare for scenarios in which renewable energy sources underperform expected output for extended periods of time. These periods of time are explicitly modeled using the crossing state stochastic model, making it a very appropriate choice for storage problems such as this one.

We saw that using the crossing state model did introduce additional complexities to even incorporate it into an SDDP algorithm (i.e. the inclusion of partially observable states) and significantly increased the length of time necessary for the solution to converge. In order to accelerate convergence, selective sampling of the outcome space was utilized to speed up the backward pass. As standard sampling tended to reduce solution quality due to its tendency to overlook high-risk, low probability outcomes, a risk-directed importance sampling algorithm was implemented to sample these elements of the outcome space with higher probability.

The resulting algorithm -- SDDP with hidden Markov models, quadratic regularization, and importance sampling in the backward pass -- converges in roughly the same time (albeit slightly longer) as the classic SDDP algorithm with the IID wind model, allowing for a direct comparison between two general approaches to developing a solution. On one hand, we can start by simplifying the problem to one that we can solve optimally and, in this case, efficiently. This is what we are doing by ignoring the complex intertemporal dynamics of the stochastic wind process and assuming an IID wind model, which lends itself to a fairly straightforward application of classic SDDP. Conversely, we can first model the problem in as much detail as possible and then use approximations in the solution algorithm to develop a practical solution. Incorporating the crossing state model and subsequently sampling to reduce CPU time does produce approximate value functions that are less optimal with respect to the unsampled version. However, the degradation in solution quality due to approximations in the solution algorithm is much less pronounced than when appoximations are made in the model that ignore key characteristics of the problem.

Finally, we also observe that augmenting the state variable to track cumulative metrics and applying an end-of-horizon threshold penalty is an effective and relatively simple way to formulate a risk-averse objective function without involving dynamic risk measures. In certain cases this is a more natural notion of risk (for example, a system can only experience so much stress over a period of time until failure) and thus a useful algorithmic tool for the practitioner. In this problem setting, the system operator is given a straightforward method for controlling the amount of acceptable cumulative shortages by adjusting the corresponding threshold parameter. Due to the success of this approach, potential future work may explore further augmentation of the state variable combined with general convex utility functions to influence solution quality.



%

\section*{Acknowledgment}

The research was supported by NSF grant CCF-1521675 and DARPA grant FA8750-17-2-0027.





\bibliographystyle{informs2014}
\bibliography{bibpaper}

\newpage

\section*{Appendix A: Details of the Crossing State Stochastic Model}

The purpose of this section is to provide sufficient details to implement the crossing state model in this problem. We will be completing the development of equations that were left for this appendix throughout the main text of the paper. This information was first presented (in greater depth) in \cite{durante2017b}, but is repeated here for the convenience of the reader. The following is presented in the context of using the model for wind power forecast errors, as is done in this paper. Though it should be noted that the model is versatile and can model other stochastic processes for which it is important to capture crossing time behavior.

Recall from Sections \ref{StateVar} and \ref{Exo Info} that we have wind power forecasts $\left\lbrace f^E_t \right\rbrace_{t=0}^T$, actual power outputs at time $t$, $E_t^W$, and forecast errors $W_t=E_t^W-f^E_t$. Assume we have access to historical training data providing us with wind power forecasts and the resulting actual output. 

First we define both up- and down-crossing times. Let the current elapsed time above forecast at time $t$, $\tau_{t}^{U}$, be
\begin{align}
\tau_{t}^{U}=\ell \text{ if }
\begin{cases}
W_{t-\ell} \leq 0\\
W_{t+\ell'} > 0 \quad\forall \ell' \in \left\lbrace 0,1,...,\ell-1 \right\rbrace.
\end{cases}
\end{align} Similarly, the current elapsed time below forecast at time $t$, $\tau_{t}^{D}$, is defined as
\begin{align}
\tau_{t}^{D}=\ell \text{ if }
\begin{cases}
W_{t-\ell} > 0 \\
W_{t+\ell'} \leq 0 \quad\forall \ell' \in \left\lbrace 0,1,...,\ell-1 \right\rbrace.
\end{cases}
\end{align} Next, let the set of all indices such that forecast errors cross over from the negative to positive regime be $\mathcal{C}^{U}=\lbrace t|W_{t-1}\leq 0 \wedge W_t >0\rbrace$. Likewise, let the set of all indices such that errors cross over from the positive to negative regime be $\mathcal{C}^{D}=\lbrace t|W_{t-1}\geq 0 \wedge W_t <0\rbrace$. The sets of all up- and down-crossings are then $\mathcal{T}^{U}=\lbrace \tau_{t}^{U}|t+1 \in \mathcal{C}^{D}\rbrace$ and $\mathcal{T}^{D}=\lbrace \tau_{t}^{D}|t+1 \in \mathcal{C}^{U}\rbrace$ respectively as a crossing time is simply a completed elapsed time above or below forecast.

For both up- and down-crossing times, there exists cumulative distribution functions $F^{U}$ and $F^{D}$ respectively. Up-crossing time distributions are quantized by partitioning into $m$ bins, splitting at the $q_i=\frac{i}{m}$ quantile points for $i=0,1,...,m-1$. An up-crossing time, $\tau_t^{U} \in \mathcal{T}^{U}$, belongs to crossing time duration bin $D_t=d$ if $q_d \leq F^{U}(\tau_t^{U}) < q_{d+1}$. Down-crossing time distributions are similarly quantized.

Our crossing state variable $I_t^{C}\equiv (C_t,D_t)$ is defined as the pair of variables describing whether or not the error is above the forecast, $C_t=\mathbf{1}_{\left\lbrace W_t>0\right\rbrace}$, and to which crossing time duration bin, $D_t$, the completed crossing time will belong. Note that this means that during the forward pass, the state $C_t$ is observable (we know if we are above or below the forecast), but the duration bin $D_t$ is not until the sign variable $C_t$ switches. However, when building the crossing state-dependent error distributions from training data for the model we can observe the duration bin at each point in time by peeking into the future to find the complete crossing time. Letting $\mathcal{I}^{C}$ be the set of all possible crossing states for the process, there exists a distribution from training data of crossing times $F^{\tau}_{i}$ for each possible crossing state $i \in \mathcal{I}^{C}$.

Transitions between crossing states are made using a transition matrix $\mathbb{P}(i'|i)$ for each pair of crossing states $(i',i) \in \mathcal{I}^{C} \times \mathcal{I}^{C}$ in which self-transitions are not allowed ($\mathbb{P}(i|i)=0$ $\forall i \in \mathcal{I}^{C}$). This matrix is computed from training data by considering only pairs of points in time $(t,t+1)$ such that $t+1 \in (\mathcal{C}^{U} \cup \mathcal{C}^{D})$ (points where the crossing state makes a transition). For all of these pairs, letting $n(I_{t+1}^{C}=i'|I_t^{C}=i)$ be the count of the transitions from state $i$ to state $i'$ occurring for each pair of crossing states $(i,i')$ and $n(I_t^{C}=i)$ be the number of times $I_t^{C}=i$ for each crossing state $i$, the empirical transition probability from crossing state $i$ to $i'$ is
\begin{align}
\label{tprob}
\mathbb{P}(i'|i)=\frac{n(I_{t+1}^{C}=i'|I_t^{C}=i)}{n(I_t^{C}=i)}.
\end{align} The duration-dependent crossing state transition probability is then a function of the elapsed time above or below forecast given by
\begin{align}
\mathbf{P}(I_{t+1}^{C}=i'|I_{t}^{C}=i,\tau_t)=
\begin{cases}
1-F_{i}^{\tau}(\tau_t) & \text{if } i'=i\\
F_{i}^{\tau}(\tau_t)\mathbb{P}(i'|i) & \text{for } i'\neq i.
\end{cases}
\end{align}

Next we describe conditioning the error generation on the crossing state. From training data, there exists empirical conditional error cumulative distribution functions $F^{W}_i$ and corresponding error density functions $\mathbf{P}(W_{t+1}|i)$ for $i \in \mathcal{I}^{C}$. Error distributions are not identical across crossing states; in fact they are likely to be quite different, such as the case where the error distribution is asymmetric. Furthermore, error distributions are likely to vary across duration bins as well. Thus, to better capture the behavior of the error process, the error generation process is conditioned on the crossing state $I_t^{C}$.

In addition to errors being crossing state-dependent, they are dependent on previous errors as well. A first order Markov chain is used to model this behavior. Similar to how the crossing time distributions are quantized, each conditional error distribution $F^{W}_i$ is partitioned into $n$ bins, splitting at the $q_j=\frac{j}{n}$ quantile points for $b=0,1,...,n^E-1$. The error $W_t$ belongs to bin $W^b_t$ if $q_b \leq F^{W}_i(W_t) < q_{b+1}$. Then, given $W_t \in W^b_t$, we form conditional empirical distributions for the error at time $t+1$ giving $\mathbf{P}(W_{t+1}|I_t^{C}, W_t^b)$.

It is important to realize that the same error $W_t$ can fall in different error bins for different crossing states. For example, the error $W_t=+50$ MW may be in bin $W_t^b=4$ for the $I_t^{C}=(1,0)$ crossing state (short up-crossings), but for the $I_t^{C}=(1,2)$ state (longer up-crossings), it may belong to bin $W_t^b=2$. To avoid additional notation, \textit{the variable $W_t^b$ is always paired with a crossing state and refers to the bin that the error $W_t$ belongs to for the corresponding crossing state}.

For each crossing state $i \in \mathcal{I}^{C}$, there also exists an error density $\mathbf{P}(W_{t+1}|i, t+1 \in \mathcal{C}^{U} \cup \mathcal{C}^{D})$. This is the distribution of the initial error given the process has just transitioned to the new crossing state $i$.

If known, the variables $I_t^{C}$, $\tau_t$, and $W^b_t$ fully determine the distribution of the exogenous information $W_{t+1}$, given by
\begin{align}
\label{Ehat}
\begin{split}
\mathbf{P}(W_{t+1}|I_t^{C}=i,\tau_t,W^b_t)=&(1-F_{i}^{\tau}(\tau_t))\mathbf{P}(W_{t+1}|i,W^b_t)
+\\
&F_{i}^{\tau}(\tau_t)\sum\limits_{i' \neq i}\mathbb{P}(i'|i)\mathbf{P}(W_{t+1}|i', t+1 \in \mathcal{C}^{U} \cup \mathcal{C}^{D}).
\end{split}
\end{align} Typically, this would imply that our post-decision information state is comprised of these variables. However, this form of the information state results in a computational issue when attempting to index VFAs on post-decision states using SDDP. Letting $\tau_i^{max}$ be the largest crossing time for crossing state $i\in \mathcal{I}^{C}$, we see that there are $\sum\limits_{i\in \mathcal{I}^{C}} n \tau_i^{max}$ possible information states at each time $t$. This number can be quite large, especially if crossing times tend to span many time periods. For this reason, we introduce the compact information state $I_t^W  \equiv \left(I_t^{C},W^b_t \right)$ which drops the current elapsed time $\tau_t$ and can only take on $2 \times m \times n$ states. This is the information state that VFAs are indexed on in the algorithm.

The corresponding pre-decision information state is given by $I_t  \equiv \left(I_t^W,E_t^W \right)$ and the function $I_t^x= S^{I,x}(I_t)$ from equation \eqref{Trans 2} is simple to define. Once $E_t^W=f_t^E+W_t$ is used to make a decision, it can be dropped from pre- to post- decision state as $W_t$ does not have a direct impact on the distribution of $W_{t+1}$. Instead we only need to maintain its error bin $W_t^b$ and the current crossing state $I_t^{C}$. Thus, we have $I_t^x=S^{I,x}(I_t)=I_t^W$. The remainder of this section is written using $I_t^x$ in all places where $I_t^W$ would also be valid. This helps connect the content in this section with the algorithms in the main text that are written for fitting VFAs to a general hidden post-decision information state $I_t^x$.

The error distributions $\mathbf{P}(W_{t+1}|i,W^b_t)$ and $\mathbf{P}(W_{t+1}|i,t+1 \in \mathcal{C}^{U} \cup \mathcal{C}^{D})$ for all $i \in \mathcal{I}^{C}$ are unaffected by the change to a compact form of $I_t^x$. However, the transition between crossing states must now be modeled with a Markov approximation of the semi-Markov model as no count of the elapsed time above or below forecast is maintained. Transition probabilities are now given by a time-invariant modified crossing state transition matrix $\tilde{\mathbb{P}}(I_{t+1}^{C}|I_t^{C})$ which allows for self-transitions. This is estimated from training data using equation \ref{tprob}; however all time periods $t$ are considered, not only pairs of points where errors switch signs. Consequently, the distribution of $W_{t+1}$, given only the compact information state $I_t^x$, is
\begin{align}
\mathbf{P}(W_{t+1}|I_t^{C}=i,W^b_t)=\tilde{\mathbb{P}}(i|i)\mathbf{P}(W_{t+1}|i,W^b_t)
+\sum\limits_{i' \neq i}\tilde{\mathbb{P}}(i'|i)\mathbf{P}(W_{t+1}|i', t+1 \in \mathcal{C}^{U} \cup \mathcal{C}^{D}).
\end{align} This is utilized in the backward pass of SDDP, for which the elapsed crossing time cannot be observed. Another important equation needed for the backward pass is
\begin{align}
\mathbf{P}(I_{t}^{C}=i|W_t)=\frac{1}{\tilde{p}^{norm}}\left(\tilde{\mathbb{P}}(i|i)\sum\limits_{b}\mathbf{P}(W_t|i,W_{t-1}^b=b)+\sum\limits_{i' \neq i}\tilde{\mathbb{P}}(i|i')\mathbf{P}(W_{t}|i, t \in \mathcal{C}^{U} \cup \mathcal{C}^{D})\right),
\end{align} where $\tilde{p}^{norm}$ is a normalization factor such that $\sum\limits_{i'\in\mathcal{I}^{C}}\mathbf{P}(I_{t}^{C}=i'|W_t)=1$. This is the probability of being in each crossing state given only an observation of $W_t$ without any elapsed time above or below the forecast. Since $I_t^x=(I_t^{C},W_t^b)$ and for each crossing state $W_t$ can only belong to one bin $W_t^b$, given we have $P(I_t^{C}=i|W_t)$ for $i \in \mathcal{I}^{C}$, we also have
\begin{align}
\mathbf{P}(I_{t}^{x}=(i,W_t^b)|W_t)=
\begin{cases}
\mathbf{P}(I_t^{C}=i|W_t) & \text{if } W_t \in W_t^b \text{ for crossing state } i\\
0 & \text{if } W_t \not\in W_t^b \text{ for crossing state } i.
\end{cases}
\end{align} for every possible post-decision information state $I_t^x$. This holds true every time we see $P(I_t^{C}=i)$ for $i \in \mathcal{I}^{C}$ in the remainder of this section as well. We can use this relation to reconcile the equations from the main text in Section \ref{SDDP Algo}, which are written for general hidden Markov models that feature a hidden post-decision information state $I_t^x$, with the equations here which are specific to the crossing state model.


In the forward pass, the current elapsed time above or below forecast, $\tau_t$, is observable. We can use this to inform our decision and thus it is placed in the knowledge state $K_t$. Given we only know the current elapsed time $\tau_t$, the future complete crossing time may belong to different crossing time duration bins $D_t$, and thus the crossing state is partially unobservable. Based on the sample path, we are able to form time $t$ beliefs about the probability the process is in each possible hidden state. This information, $\left\lbrace\mathbf{P}(I_t^{C}=i)\right\rbrace_{i\in\mathcal{I}^{C}}$, along with $\tau_t$ and $W_t$, forms the knowledge state $K_t$.

Given $K_t$, we can derive our belief about the distribution of the error at time $t+1$, $\mathbf{P}(W_{t+1}|K_t)$. We are able to determine the sign of the error, $C_t$, based on $W_t$. Then, for each possible value of $D_t$, and corresponding crossing state $I_t^{C}=(C_t,D_t)$, $W_t$ can belong to only one error bin $W_t^b$. We then have
\begin{align}
\label{Etplus1givenKt}
\mathbf{P}(W_{t+1}|K_t)=\sum\limits_{i\in \mathcal{I}^{C}} \mathbf{P}(I_t^{C}=i)\mathbf{P}(W_{t+1}|i,\tau_t,W^b_t),
\end{align}
where $\mathbf{P}(W_{t+1}|i,\tau_t,W^b_t)$ is given by equation \ref{Ehat}.

Subsequently, following the observation of $W_{t+1}$, a Bayesian update is performed on the knowldge state at each time step according to the update function $(K_{t+1},I_{t+1})=S^{K,W}(K_t,I_t^x,W_{t+1})$ (equation \eqref{Trans 4}) defined by the following two cases:

\begin{hitemize}
\item Case 1: $sign(W_{t+1})=sign(W_{t})$. In this case $\tau_{t+1}=\tau_{t}+1$. This is then used to compute the likelihood that the future completed crossing time belongs to bin $D_t$ for each crossing state $i=(C_t,D_t)$ given the elapsed time above or below forecast $\tau_{t+1}$. This likelihood is given by $1-F^{\tau}_{i}(\tau_{t+1})$. Furthermore, the likelihood of observing error $W_{t+1}$ in crossing state $i$ given a crossing state transition has not occurred is $\mathbf{P}(W_{t+1}|i,W^b_t)$. Thus with prior beliefs, $\mathbf{P}(I_{t}^{C}=i)$ for $i\in\mathcal{I}^{C}$, we compute the posterior beliefs
\begin{equation}
\mathbf{P}(I_{t+1}^{C}=i)=\frac{1}{p^{norm}}\left(\mathbf{P}(I_{t}^{C}=i)(1-F^{\tau}_{i}(\tau_{t+1}))\mathbf{P}(W_{t+1}|i,W^b_t)\right),
\end{equation} where $p^{norm}=\sum\limits_{i'\in\mathcal{I}^{C}}\mathbf{P}(I_{t}^{C}=i')(1-F^{\tau}_{i'}(\tau_{t+1}))\mathbf{P}(W_{t+1}|i',W^b_t)$.

\item Case 2: $sign(W_{t+1})\neq sign(W_t)$. In this case $\tau_{t+1}=1$ and we are able to determine the crossing state at time $t$ based on the sign of $W_t$ and the completed crossing time $\tau_t$; let this be state $i^*$. Additionally, we know that a crossing state transition has taken place. The likelihood of observing error $W_{t+1}$ in crossing state $i$ given a crossing state transition has just occurred is $\mathbf{P}(W_{t+1}|i, t+1 \in \mathcal{C}^{U} \cup \mathcal{C}^{D})$. Thus, for $i\in\mathcal{I}^{C}$, posterior beliefs are given by
\begin{equation*}
\mathbf{P}(I_{t+1}^{C}=i)=\frac{1}{p^{norm}}\left(\mathbb{P}(i|i^*)\mathbf{P}(W_{t+1}|i, t+1 \in \mathcal{C}^{U} \cup \mathcal{C}^{D})\right),
\end{equation*} where $\mathbb{P}(i|i^*)$ is defined in equation \eqref{tprob} and $p^{norm}=\sum\limits_{i'\in\mathcal{I}^{C}}\mathbb{P}(i'|i^*)\mathbf{P}(W_{t+1}|i',t+1 \in \mathcal{C}^{U} \cup \mathcal{C}^{D})$.
\end{hitemize} In either case, to complete the transition function we also add $E_{t+1}^W=f_{t+1}^E+W_{t+1}$ back to the pre-decision information state $I_{t+1}$. Given these recursive updating formulas for the knowledge state, we only need to initialize our beliefs at $t=0$. Given $W_0$, we use a discrete uniform distribution for the initial beliefs: $\mathbf{P}(I_0^{C}=i)=1/m$ for $i \in \mathcal{I}^{C}$ such that $C_0=sign(W_0)$. Setting $\tau_0=1$, this forms $K_0$.

\end{document}